\documentclass[preprint,3p,12pt]{elsarticle}
\biboptions{numbers,sort&compress}

\usepackage{amssymb,verbatim,amsmath}
\usepackage{hyperref}
\usepackage{graphicx}
\usepackage{caption}
\usepackage{multirow}
\usepackage{pdflscape}
\usepackage{float}

\usepackage{xcolor}

\usepackage{longtable}

\newtheorem{remark}{Remark}

\numberwithin{equation}{section}
\numberwithin{remark}{section}

\def\dhp#1{\mathop {#1}\limits_{+h}}
\def\dtp#1{\mathop {#1}\limits_{+\tau}}
\def\dtm#1{\mathop {#1}\limits_{-\tau}}
\def\dsp#1{\mathop {#1}\limits_{+s}}
\def\dsm#1{\mathop {#1}\limits_{-s}}
\def\dhm#1{ \mathop{#1}\limits_{-h}}

\newcommand{\energconservtxt}{\text{Energy conservation}}
\newcommand{\mmomconservtxt}{Mass and momentum conservation;}
\newcommand{\empconservtxt}{\bgroup\def\arraystretch{1.25}{\begin{tabular}[c]{@{}l@{}}Energy, momentum\\and mass conservation;\end{tabular}}\egroup}
\newcommand{\noaddcls}{\bgroup\def\arraystretch{1.25}{\begin{tabular}[c]{@{}l@{}}Energy, momentum\\and mass conservation;\\No additional conservation\\laws\end{tabular}}\egroup}
\newcommand{\noaddclsx}{\bgroup\def\arraystretch{1.25}{\begin{tabular}[c]{@{}l@{}}Mass and energy conservation;\\No additional conservation laws\end{tabular}}\egroup}

\newcommand{\energyOrMomentumText}{\bgroup\def\arraystretch{1.25}\begin{tabular}[c]{@{}l@{}}
Mass and momentum conservation\\
\qquad\emph{or}\\
mass and energy conservation ---\\
see schemes in~\cite{dorodnitsyn2019shallow}
\end{tabular}
\egroup}

\def\arraystretch{1.75}

\definecolor{light-gray}{gray}{0.5}

\begin{document}

\title{Invariant conservative difference schemes
for shallow water equations
in Eulerian and Lagrangian coordinates}

\begin{frontmatter}
\author[mymainaddress]{V.~A. Dorodnitsyn}
\ead{Dorodnitsyn@keldysh.ru,dorod2007@gmail.com}

\author[mymainaddress,mysecondaryaddress]{E.~I. Kaptsov\corref{mycorrespondingauthor}}
\cortext[mycorrespondingauthor]{Corresponding author}
\ead{evgkaptsov@gmail.com}

\address[mymainaddress]{Keldysh Institute of Applied Mathematics,\\
Russian Academy of Science, Miusskaya Pl. 4, Moscow, 125047, Russia}

\address[mysecondaryaddress]{School of Mathematics, Institute of Science, \\
Suranaree University of Technology, 30000, Thailand}


\begin{abstract}
The one-dimensional shallow water equations in Eulerian coordinates
are considered. Relations between symmetries and conservation laws
for the potential form of the equations, and symmetries and conservation laws
in Eulerian coordinates are shown. An invariant difference scheme
for equations in Eulerian coordinates with arbitrary bottom topography
is constructed. It possesses all the finite-difference analogues of the conservation
laws. Some bottom topographies require moving meshes in Eulerian coordinates,
which are stationary meshes in mass Lagrangian coordinates. The developed
invariant conservative difference schemes are verified numerically
using examples of flow with various bottom topographies.
\end{abstract}

\begin{keyword}
Shallow water equations\sep Eulerian and Lagrangian coordinates \sep
Lie point symmetries \sep conservation law \sep invariant numerical
scheme
\end{keyword}

\end{frontmatter}

\medskip

\begin{center}
\today
\end{center}

\medskip

\section{Introduction}

The one-dimensional shallow water equations in Eulerian and
Lagrangian coordinates are considered. The recent
article~\cite{dorodnitsyn2019shallow} was devoted to the group
symmetries and conservation laws preservation in finite-difference
modeling for the shallow water equations in  Lagrangian coordinates.
The present article is devoted to the construction of invariant
conservative difference schemes for the shallow water equations in
Eulerian coordinates.

The shallow water equations describe the flow below a free surface in a fluid.
They are commonly used to model processes in water basins, atmosphere, tidal oscillations and
gravity waves~\cite{bk:Whitham[1974],bk:Ovsyannikov[2003],bk:Vallis[2006],bk:PetrosyanBook[2010]}.
In particular, one-dimensional shallow water equations
are widely used in modeling of transient open-channel flow and
surface runoff.

The construction of nontrivial exact solutions of the shallow water
equations is a rather difficult problem even in the one-dimensional
case~(some exact solutions can be found in~\cite{bk:Camassa2019,bk:PetrosyanBook[2010],bk:Bernetti[2008],bk:HanHantke[2012]}).
The nonlinearity of the equations and the absence of their exact solutions
emphasize the importance of numerical
modeling of the shallow water equations, which has been the subject
of many publications,
e.g.~\cite{bk:YeleninKrylov[1982], bk:Bihlo_numeric[2012],
bk:Bihlo_numeric[2017],bk:Bihlo_numeric[2019],
bk:MurshedFutai[2019],bk:KhakimzyanovIV,
bk:DyakonovaKhoperskov,bk:MoralesCastro,dorodnitsyn2019shallow}.

Depending on the problem studied, the shallow water equations are
considered as in Lagrangian coordinates associated with the movement
of particles of the medium, or in Eulerian coordinates associated
with a fixed bottom. In~\cite{dorodnitsyn2019shallow}, the authors considered the
equations in Lagrangian coordinates. The present paper is
devoted to the construction of invariant conservative difference
schemes for the shallow water equations in Eulerian coordinates.

Group analysis, based on fundamental theory by Sophus
Lie~\cite{Lie15, bk:Lie[1891b], bk:Lie-Scheffers[1896]}, has proven to be
an effective tool to study equations of mathematical physics and
continuous mechanics, to investigate their qualitative properties,
conservation laws and exact
solutions~\cite{bk:Ovsyannikov[1962],bk:Olver,bk:Ibragimov1985,bk:Bluman1989,
bk:HandbookLie_v1,bk:Gaeta1994}. In particular, group analysis and
the group classifications of the shallow water equations can be
found
in~\cite{bk:HandbookLie_v2,bk:LeviNicciRogersWint[1989],bk:ClarksonBila[2006],
AksenovDruzkov2016,bk:AksenovDruzkov_classif[2019],
bk:KaptsovMeleshko_1D_classf[2018],
bk:Andronikos2019,
bk:MeleshkoSW2020,bk:MeleshkoSamatova2020}.
Some extended nonlinear models (such as the Green--Naghdi equations)
were considered in~\cite{bk:SiriwatKaewmaneeMeleshko2016,
bk:SzatmariBihlo[2014]}.

It turned out that many tools for group analysis of differential
equations can be applied to finite-difference
equations~\cite{Maeda1, Maeda2,Dor_1, Dor_2, Dor_3,
bk:DorodKozlovWint[2004],[LW-2],bk:DorodKozlovWinternitz[2000],
Quisp, [LW-3], bk:Dorodnitsyn[2011], Vinet,
bk:Hydon_book[2014],bk:DorodKozlovWintKaptsov[2015]}. An important
peculiarity of group analysis in finite difference spaces is the
need to take into account the nonlocality of difference operators
and the geometry of difference
meshes~\cite{Dor_1,bk:Dorodnitsyn[2011]}. Invariant
finite-difference schemes, i.e. difference equation and a mesh,
being constructed by means of difference invariants, admit the same
symmetries as their differential counterparts. They can also possess
difference analogues of conservation laws and difference invariant
solutions~\cite{DORODNITSYN2019201,KOZLOV2019,bk:Kozlov[2007],dorodnitsyn2019shallow}.
Invariant difference schemes of ordinary differential
equations can be also reduced and completely integrated in
certain cases~\cite{bk:Dorodnitsyn[2011],bk:DorodKaptsov[2013],bk:DorodKozlovWintKaptsov[2014]}.
The finite-difference analogues of Lagrangian
and the Hamiltonian formalism were developed in~\cite{Dor_3,bk:DorodKozlovWint[2004],bk:Dorod_Hamilt[2011]}
and~\cite{bk:Dorod_Hamilt[2010],bk:Dorod_Hamilt[2011]}.
Moreover, absent the Lagrangian and Hamiltonian,
there is a more general approach based on the Lagrange operator identity
and adjoint equation method~\cite{[Bluman1]}, the difference analogue of which was also developed in \cite{bk:DorodKozlovWintKaptsov[2014],bk:DorodKozlovWintKaptsov[2015]}.
An alternative ``direct method'' of local conservation law computation~\cite{bk:BlumanAnco_adjoint[1996],bk:BlumanCheviakovAnco}
which employs Euler differential operators was adapted
and applied to the case of finite-difference equations in~\cite{dorodnitsyn2019shallow,bk:ChevDorKap2020,bk:DorKapMelGN2020}.

\medskip

In the present paper the relationship between symmetries and conservation laws for the potential form of the equations and symmetries and conservation laws in Eulerian variables is studied. For equations in Eulerian coordinates with arbitrary bottom topography, an invariant difference scheme is constructed. For some special bottom shapes, moving meshes in Eulerian coordinates which are stationary meshes in mass Lagrangian coordinates are needed. The developed invariant conservative difference schemes are tested numerically on examples
of flow over various  bottom topography.

\medskip

The paper is organized as follows.
In Section 2 the shallow water equation
in Eulerian coordinates are considered. Application of the
Noether theorem to the potential form of the shallow water equations
yields conservation laws which can be transformed back to physical
variables in Eulerian coordinates.
It Section 3, invariance of differential and finite-difference shallow water
equations considered for various special cases of bottom shapes when
the admitted symmetry group has an extension is presented. For some cases the simple
uniform orthogonal mesh is not applicable, and all spacial cases are
collected in Table~\ref{tab:HtableDelta}.
Invariant conservative difference schemes
for an arbitrary bottom shape are constructed in Section~4.
In Sections~5 and~6 the numerical implementation of
the obtained invariant schemes is
performed.
The results are summarized in Conclusion.

\section{Shallow water equation for Eulerian coordinates}

\begin{figure}[ht]
\centering
\includegraphics[width=0.25\linewidth]{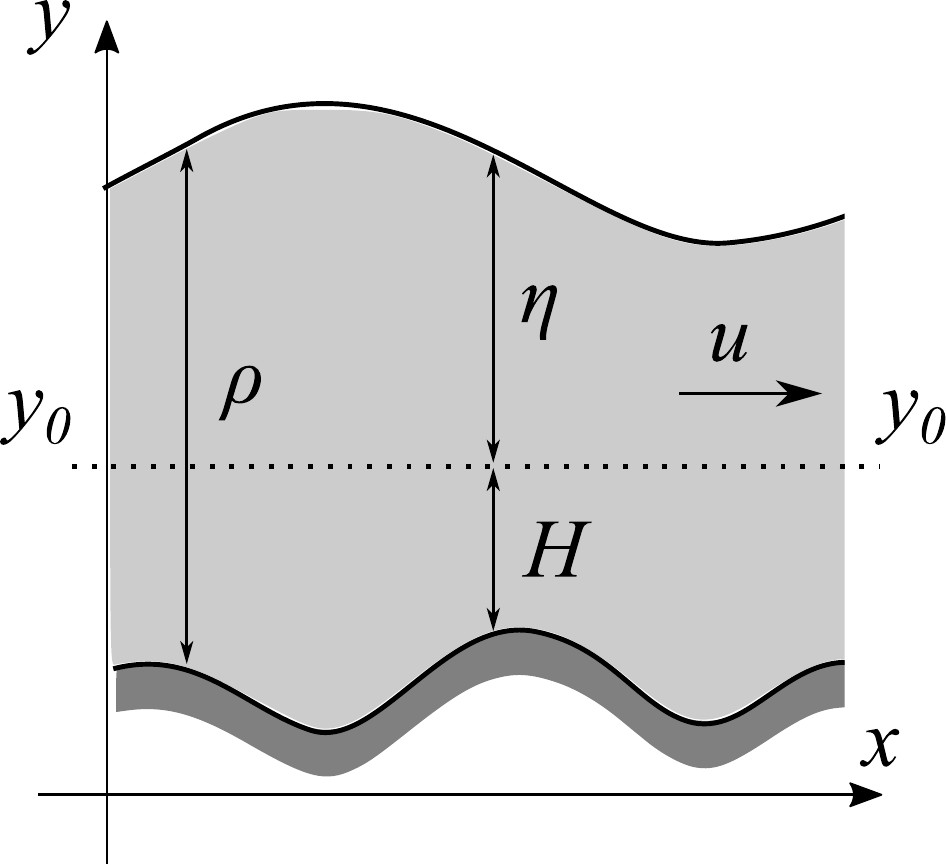}
 \caption{The one-dimensional shallow water flow over uneven bottom.
 The variable~$\rho = \eta + H$ which is introduced below
 describes the depth of the fluid.}
\label{fig:notation}
\end{figure}

The system of the one-dimensional shallow water equations with an uneven bottom
topography in Eulerian coordinates has the following form
\begin{equation} \label{Euler1}
\eta_t + ((\eta +H)u)_x = 0,
\end{equation}
\begin{equation} \label{Euler2}
u_t +uu_x + \eta_x=0,
\end{equation}
where $u(t,x)$ is the velocity of the particles in continuous medium,
$\eta(t,x)$ is the height of the fluid over the chosen undisturbed level~$y_0$~(see Figure~\ref{fig:notation}),
and the bottom topography is described by the function~$H=H(x)$.

\bigskip

One can reduce the linear bottom case~($H(x) = k x$, where $k \neq 0$ is constant)
to the flat bottom~($H = 0$) by the following change
of variables~\cite{bk:ChirkunovPikmullina[2014]}
\begin{equation}
\label{bottom_tr_Chirkunov}
    t = \frac{t_*}{k},
    \quad
    x = \frac{1}{k} \left( x_* + \frac{t_*^2}{2}\right),
    \quad
    u = u_* + t_*,
    \quad
    \eta = \eta_* - x_* - \frac{t_*^2}{2}.
\end{equation}
In case of the flat bottom, it turns out that
system~(\ref{Euler1}), (\ref{Euler2}) possesses a especially simple form.
The hodograph transformation
\begin{equation}\label{hodograph_Euler}
\begin{array}{c}
    x = x(\eta, u),
    \qquad
    x_\eta = -\frac{u_t}{\Delta},
    \qquad
    x_u = \frac{\eta_t}{\Delta},
    \\
    t = t(\eta, u),
    \qquad
    t_\eta = \frac{u_x}{\Delta},
    \qquad
    t_u = -\frac{\eta_x}{\Delta},
    \\
    \Delta = \eta_t u_x - \eta_x u_t \neq 0
\end{array}
\end{equation}
allows one to linearize~\cite{bk:YanenkRojd[1968],bk:Ovsyannikov[2003]}
the system and arrive at the following equations
\[
    x_u - u t_u + {\eta} t_{\eta} = 0,
\]
\[
    x_{\eta} + t_u - u t_{\eta} = 0.
\]

The latter system is linear and has infinite number
of symmetries, which means that system~(\ref{Euler1}),~(\ref{Euler2}) also has this property.

\bigskip


\medskip
\subsection{Application of Noether's theorem}

In order to find conservation laws of
system~(\ref{Euler1}),~(\ref{Euler2}), one can introduce a potential
function and derive a potential form of the equations. The
potential form allows one to represent the system as Euler--Lagrange
equations for some Lagrangian, and then to obtain
conservation laws by means of Noether's theorem.

One introduces the potential for equatios~(\ref{Euler1}),~(\ref{Euler2})
\begin{equation} \label{Pot}
u=w_x.
\end{equation}
Then one has that
\begin{equation} \label{PotentialE}
w_{tx}+w_xw_{xx}+\eta_x=0,
\end{equation}
or
\begin{equation} \label{PE}
D_x[w_t+\frac{w_x^2}{2}+\eta]=0.
\end{equation}

An integration of the latte equation yields
\begin{equation} \label{ialE}
[w_t+\frac{w_x^2}{2}+\eta]=K(t),
\end{equation}
where $K(t)$ is an arbitrary function, which can be considered as
zero (by means of a change of variables).
Equation~ (\ref{ialE}) with $K=0$ becomes
\begin{equation} \label{23}
\eta = -w_t - \frac{w_x^2}{2}.
\end{equation}
Substituting $\eta$ into equation~(\ref{Euler1}), one arrives at
the equation
\begin{equation} \label{233}
-w_{tt}-2w_x w_{tx}+ w_{xx}\left(H(x)-w_t-\frac{3}{2}w_x^2
\right)+H_x w_x=0
\end{equation}
which possesses the following symmetries
\begin{equation} \label{symmetries}
 X_1=\frac{\partial}{\partial t}, \qquad
 X_w=\frac{\partial}{\partial w}.
\end{equation}
Equation (\ref{233}) is an Euler-Lagrange equation with the Lagrangian
\begin{equation} \label{234}
{\cal{L}} = \frac{w_t^2}{2} -
H(x)\frac{w_x^2}{2}+\frac{w_x^4}{8}+w_t \frac{w_x^2}{2}.
\end{equation}
Now one can apply the Noether theorem to obtain conservation laws.

\bigskip

Using the translation symmetry
\begin{equation} \label{5234}
X_1=\frac{\partial}{\partial t}
\end{equation}
one finds the conservation law
\begin{equation} \label{2341}
D_t\left[\frac{w_t^2}{2}-H(x)\frac{w_x^2}{2}+ \frac{w_x^4}{8}\right] -
D_x\left[-w_t\left(-H(x)w_x + \frac{w_x^3}{2}+ w_tw_x\right)\right]=
\end{equation}
$$
 =w_t(-w_{tt}-2w_x
w_{tx}+ w_{xx}\left(H(x)-w_t-\frac{3}{2}w_x^2 \right)+H_xw_x)=0.
$$

\bigskip
Returning to Eulerian coordinates by means of the substitution
$$
w_x =u,  \quad  w_t=-\eta  - \frac{u^2}{2},
$$
one has
$$
D_t\left[-\frac{1}{2}\left(\eta +\frac{u^2}{2}\right)^2-H(x)\frac{u^2}{2}+
\frac{u^4}{8} \right] +
$$
$$
+ D_x\left[\left(\eta +\frac{u^2}{2}\right)\left(-H(x)u +
\frac{u^3}{2}-u\left(\eta +\frac{u^2}{2}\right)\right)\right]=
$$
\begin{equation} \label{Eul1}
-\frac{1}{2}\{D_t\left[u^2(\eta +H)+ \eta^2\right]
+ D_x\left[u(\eta +H)(u^2 +2\eta)\right]\}=
\end{equation}
$$
=- \left(\eta +\frac{u^2}{2}\right)(\eta_t +[(\eta+H)u]_x) - u(\eta + H)
(u_t+uu_x+\eta_x)=0.
$$
The latter conservation law exists for any $H(x)$ and coincides
with the conservation law found by direct method applied
to system~(\ref{Euler1}),~(\ref{Euler2})~\cite{bk:AksenovDruzkov_classif[2019]}.

\bigskip

Applying Noether's theorem to the symmetry
\begin{equation} \label{5rt4}
X_w=\frac{\partial}{\partial w}
\end{equation}
one derives the conservation law
\begin{equation} \label{2CL}
D_t[w_t+\frac{w_x^2}{2}] + D_x[-H(x)w_x + \frac{w_x^3}{2}+ w_tw_x]=
\end{equation}
$$
 =-\left(-w_{tt}-2w_x
w_{tx}+ w_{xx}\left(H(x)-w_t-\frac{3}{2}w_x^2 \right)+H_xw_x\right)=0.
$$

In Eulerian coordinates this gives the original equation:
equations:
\begin{equation} \label{euler2CL}
D_t[-\eta] + D_x\left[-H(x)u + \frac{u^3}{2}- u\left(\eta + \frac{u^2}{2}\right)\right]
=-(\eta_t + ((\eta +H)u)_x) = 0.
\end{equation}

For some types of a bottom shapes, system~(\ref{Euler1}), (\ref{Euler2})
possesses additional symmetries for which Noether's theorem cannot be applied.
We will discuss it later.

\section{Invariance of differential and finite-difference shallow water equations}

In general, the Lie group of transformations admitted by
system~(\ref{Euler1}), (\ref{Euler2}) is highly depend on the
particular form of the bottom~$H$.
In case of an arbitrary bottom topography~$H=H(x)$,
equations~(\ref{Euler1}),(\ref{Euler2}) admit the
single generator
\begin{equation}
  X_1 = \frac{\partial}{\partial t},
\end{equation}
which composes the kernel of the admitted Lie algebras.

The group classification of one-dimensional shallow water equations
over uneven bottom in Eulerian coordinates was performed in~\cite{AksenovDruzkov2016}
(one can also find the corresponding
group classifications in Lagrangian coordinates
in~\cite{bk:AksenovDruzkov_classif[2019],bk:KaptsovMeleshko_1D_classf[2018]}).
Further we consider these results in a more general form taking
into account the equivalence transformations~\cite{bk:Ovsyannikov[1962]}
which are obtained below.

We seek for the group of equivalence transformations by considering
the generator
\[
X^\textrm{e} = \xi^t \frac{\partial}{\partial t}
    + \xi^x \frac{\partial}{\partial x}
    + \zeta^u \frac{\partial}{\partial u}
    + \zeta^\eta \frac{\partial}{\partial \eta}
    + \zeta^H \frac{\partial}{\partial H},
\]
where $\xi^t$, $\xi^x$, $\zeta^u$ and $\zeta^\eta$ are functions of $t,x$ and $u$,
and $\zeta^H$ depends on $t,x,u$ and~$H$.
Performing the standard procedure~\cite{bk:Ovsyannikov[1962],bk:HandbookLie_v2}, one gets the
following group of equivalence transformations of system~(\ref{Euler1}),~(\ref{Euler2})
\begin{equation} \label{eqtrans}
  \begin{array}{c}
    X^\textrm{e}_1 = \frac{\partial}{\partial t},
    \qquad
    X^\textrm{e}_2 = \frac{\partial}{\partial x},
    \qquad
    X^\textrm{e}_3 = \frac{\partial}{\partial \eta} - \frac{\partial}{\partial H},
    \\
    X^\textrm{e}_4 =  t \frac{\partial}{\partial t} + x \frac{\partial}{\partial x},
    \qquad
    X^\textrm{e}_5 =
        x \frac{\partial}{\partial x}
        + u \frac{\partial}{\partial u}
        + 2\eta \frac{\partial}{\partial \eta}
        + 2 H \frac{\partial}{\partial H}.
  \end{array}
\end{equation}
One can easily check the involutions
\begin{equation} \label{eqtransD}
  x \mapsto -x, \quad t \mapsto -t,
  \qquad
  \text{and}
  \qquad
  x \mapsto -x, \quad u \mapsto -u
\end{equation}
are also admitted by the system.

According to the group classification~\cite{AksenovDruzkov2016},
the list of bottoms and the corresponding symmetries admitted
by the shallow water equations up to the actions
of the group of equivalence transformations~(\ref{eqtrans}),~(\ref{eqtransD}),
are provided in Table~\ref{tab:Htable}.
In Table~\ref{tab:HtableDelta} the admitted symmetries and conservation laws
for the finite-difference schemes in Eulerian and Lagrangian coordinates
are given.

\begin{figure}[ht]
\centering
\includegraphics[width=0.25\linewidth]{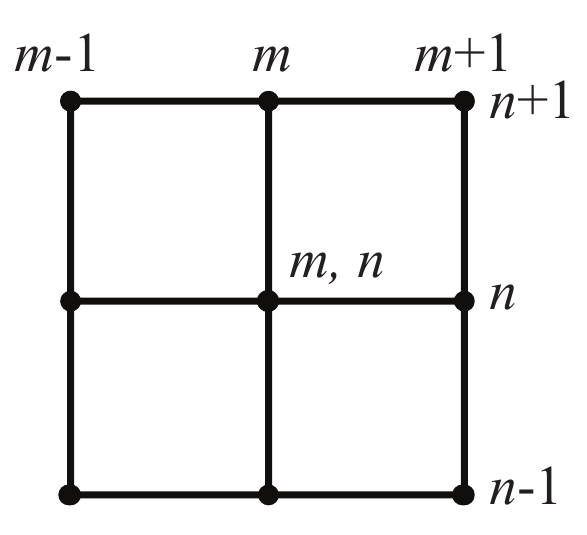}
 \caption{9-point stencil.}
\label{fig:9pt-template}
\end{figure}

We seek for invariant finite-difference equations defined on
simple orthogonal regular meshes. For such equations it is necessary
to preserve the following orthogonality and uniformness conditions
(see~\cite{bk:Dorodnitsyn[2011]}) under the actions of the
generators
\begin{equation} \label{mesh_conds_uni}
  \dhp{D}\dhm{D}(\xi^x) = 0,
  \qquad
  \dtp{D}\dtm{D}(\xi^t) = 0,
\end{equation}
\begin{equation} \label{mesh_conds_ortho}
  \dhp{D}(\xi^t) = -\dtp{D}(\xi^x),
\end{equation}
where $\underset{\pm\tau}{D}$ and $\underset{\pm{h}}{D}$ are
finite-difference differentiation operators
\[
    \underset{+\tau}{D} = \frac{\underset{+\tau}{S} - 1}{t_{n+1} - t_{n}},
    \quad
    \underset{-\tau}{D} = \frac{1 - \underset{-\tau}{S}}{t_n - t_{n-1}},
    \quad
    \underset{+h}{D} = \frac{\underset{+h}{S} - 1}{x_{m+1} - x_{m}},
    \quad
    \underset{-h}{D} = \frac{1 - \underset{-h}{S}}{x_m - x_{m-1}},
\]
which are defined through the finite-difference shifts
\[
    \def\arraystretch{1.5}
    \begin{array}{c}
    \displaystyle
    \underset{\pm\tau}{S}(f(t_n, x_m, u^{n}_{m}, \eta^{n}_{m}))
        = f(t_{n \pm 1}, x_m, u^{n \pm 1}_{m}, \eta^{n \pm 1}_{m}),
    \\
    \displaystyle
    \underset{\pm{h}}{S}(f(t_n, x_m, u^{n}_{m}, \eta^{n}_{m}))
        = f(t_n, x_{m\pm 1}, u^{n}_{m \pm 1}, \eta^{n}_{m \pm 1}).
    \end{array}
\]
The indices $n$ and $m$ are changed along time and space
axes~$t$ and~$x$ correspondingly~(further we consider a 9-point
stencil which is depicted in Figure~\ref{fig:9pt-template}).
For brevity, the following notation~\cite{bk:SamarskyPopov_book[1992]} is also used
for the variables
\begin{equation}\label{SamarskiyNotation}
\begin{array}{c}
    \big(
        t_n, t_{n-1}, t_{n+1};
        x_m, x_{m-1}, x_{m+1};
        u^n_m, u^n_{m-1}, u^n_{m+1},
        u^{n-1}_m, u^{n-1}_{m-1}, u^{n-1}_{m+1},
        u^{n+1}_m, u^{n+1}_{m-1}, u^{n+1}_{m+1};
        \\
        \eta^n_m, \eta^n_{m-1}, \eta^n_{m+1},
        \eta^{n-1}_m, \eta^{n-1}_{m-1}, \eta^{n-1}_{m+1},
        \eta^{n+1}_m, \eta^{n+1}_{m-1}, \eta^{n+1}_{m+1}
    \big)
    \\
    \equiv
    (
        t, \check{t}, \hat{t};
        x, x_-, x_+;
        u, u_+, u_-,
        \check{u}, \check{u}_-, \check{u}_+;
        \hat{u}, \hat{u}_-, \hat{u}_+,
        \eta, \eta_-, \eta_+,
        \check{\eta}, \check{\eta}_-, \check{\eta}_+
        \hat{\eta}, \hat{\eta}_-, \hat{\eta}_+,
    )
\end{array}
\end{equation}
and for their derivatives
\[
    \displaystyle
    u_t = \underset{+\tau}{D}(u),
    \qquad
    u_{\check{t}} = \underset{-\tau}{D}(u),
    \qquad
    u_{t\check{t}} = \underset{+\tau}{D}\underset{-\tau}{D}(u),
    \qquad
    u_x = \underset{+x}{D}(u),
    \qquad
    u_{\bar{x}} = \underset{-h}{D}(u),
\]
\[
    \displaystyle
    u_{x\bar{x}} = \underset{+h}{D}\underset{-h}{D}(u),
    \qquad
    \eta_t = \underset{+\tau}{D}(\eta),
    \qquad
    \eta_{\check{t}} = \underset{-\tau}{D}(\eta),
    \qquad
    \eta_{t\check{t}} = \underset{+\tau}{D}\underset{-\tau}{D}(\eta),
    \qquad
    \text{etc.}
\]

\medskip

Most of the generators in Eulerian coordinates do not satisfy conditions~(\ref{mesh_conds_ortho}).
As an example, consider the generators $X^1_4$ and $X^1_5$ (see Table~\ref{tab:HtableDelta}):
\begin{equation}
\begin{array}{c}
  X^1_4:
  \quad
  \dhp{D}(\xi^t) = \dhp{D}(0) = 0,
  \quad
  \dtp{D}(\xi^x) = \dtp{D}(t) = 1,
  \\
  \dhp{D}(\xi^t) \neq -\dtp{D}(\xi^x),
  \quad \textrm{i.~e.,} \quad
  0 \neq -1,
\end{array}
\end{equation}
and the generator $X^1_4$ breaks orthogonality of the mesh.
The same is true for the generator $X^1_5$:
\begin{equation}
\begin{array}{c}
  X^1_5:
  \quad
  \dhp{D}(\xi^t) = \dhp{D}\left( \frac{1}{2}x - \frac{3}{2} tu \right)
  \neq
  -\dtp{D}(\xi^x) = -\dtp{D}\left( \frac{3}{2}t\eta - \frac{3}{4}tu^2 \right).
\end{array}
\end{equation}
(Actually, the orthogonality can hold for particular solutions $u$, $\eta$, but
in general it evidently does not hold).
In a similar way one can verify that the generators
$X^3_3$, $X^3_4$,
$X^4_3$, $X^4_4$,
and, in general,
$X^1_\infty$ break the mesh orthogonality as well.

Notice that for the arbitatry bottom case and for the cases \#\#5--7 (see Table~\ref{tab:HtableDelta}),
one can construct invariant schemes
on the uniform orthogonal mesh
which admit the same Lie
algebra as the original differential system.
Thus, it remains to consider the cases of flat, linear and parabolic bottoms.
Actually, the linear bottom case can be reduced to the case of a flat bottom by
the change of variables~(\ref{bottom_tr_Chirkunov}), so we further
restrict our consideration to the flat and parabolic bottoms.

\begin{figure}[ht]
\begin{center}
  \begin{minipage}[c]{0.3\linewidth}
    \includegraphics[width=\linewidth]{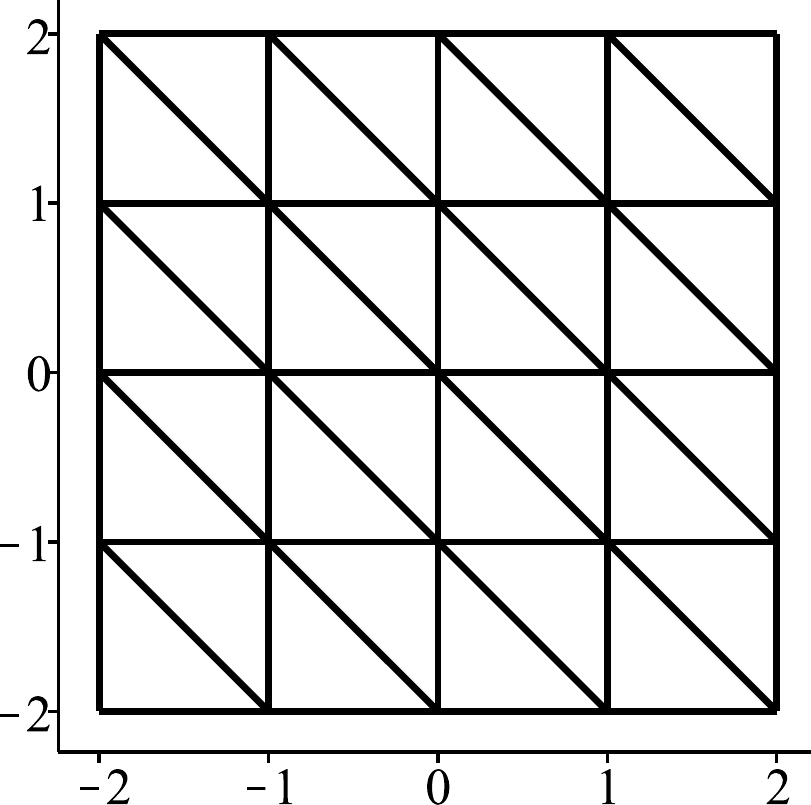}
  \end{minipage}
  \hfil
  \begin{minipage}[c]{0.3\linewidth}
    \includegraphics[width=\linewidth]{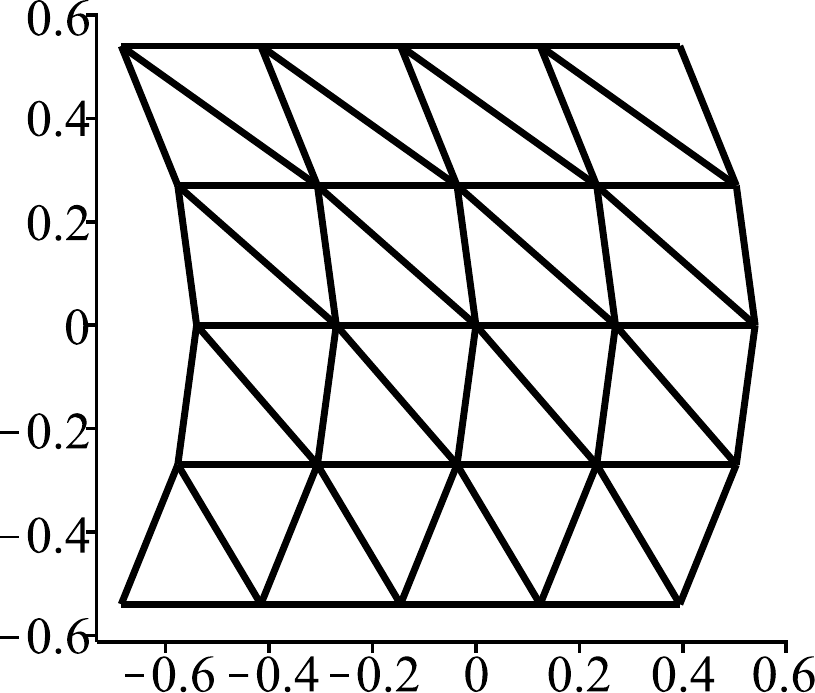}
  \end{minipage}
\caption{Transformation~(\ref{bottom_tr_Chirkunov})
(for $k=0.27$) affects orthogonality of the mesh (in the $(x,t)$-space).\\
Left:~the original mesh. Right:~the transformed mesh.}
\label{fig:transform1}
\end{center}
\end{figure}

\begin{remark}
One can check that the scalar product $\vec{h} \cdot \vec{\tau} = 0$
is not preserved by point transformation~(\ref{bottom_tr_Chirkunov}):
\[
    \vec{h}_* \cdot \vec{\tau}_* = -\frac{k^3}{2} h \tau (2 t + \tau) \neq 0,
\]
i.e., it does not preserve orthogonality of a mesh (see in Figure~\ref{fig:transform1}).
This is the reason why in the finite-difference space the group
admitted in the flat bottom case does not coincide with that for
the case of a linear bottom.
\end{remark}

\subsection{Flat bottom case}

Consider the following 5-parametric algebra which is admitted by the shallow water equations with flat bottom
\begin{equation}
  \begin{array}{c}
    X_1 = \frac{\partial}{\partial t},
    \qquad
    X_2 = \frac{\partial}{\partial x},
    \qquad
    X_3 = t \frac{\partial}{\partial t} + x \frac{\partial}{\partial x},
    \\
    X_4 = x \frac{\partial}{\partial x} + u \frac{\partial}{\partial u} + 2\eta \frac{\partial}{\partial \eta},
    \qquad
    X_5 = t \frac{\partial}{\partial x} + \frac{\partial}{\partial u}.
 \end{array}
\end{equation}
One can see that the generator $X_5$ does not satisfy criterion~(\ref{mesh_conds_ortho}),
i.e., the orthogonal mesh is not invariant under transformation with this generator.
Thus, one should look for an invariant moving mesh in Eulerian coordinates.
The finite-difference invariants of the generators~$X_1$--$X_5$
in 14-dimensional space
\begin{equation}\label{14stencil}
  (t,\hat{t}, x, x_+ \hat{x}, \hat{x}_+, u, u_+ \hat{u}, \hat{u}_+, \eta, \eta_+ \hat{\eta}, \hat{\eta}_+)
\end{equation}
are
\begin{equation} \label{flatInvs}
\begin{array}{c}
  {\displaystyle
  \frac{\hat{x}_+ - \hat{x}}{h_+},
  \qquad
  \frac{1}{\sqrt{\eta}} \frac{h_+}{\tau},
  \qquad
  \frac{1}{\sqrt{\eta}}\left( \frac{\hat{x} - x}{\tau} - u \right),}
  \\
  {\displaystyle
  \frac{1}{\sqrt{\eta}}(u_+ - u),
  \qquad
  \frac{1}{\sqrt{\eta}}(\hat{u} - u),
  \qquad
  \frac{1}{\sqrt{\eta}}(\hat{u}_+ - \hat{u}),}
  \\
  {\displaystyle
  \frac{\hat{\eta}}{\eta},
  \qquad
  \frac{\hat{\eta}_+}{\eta},
  \qquad
  \frac{\eta_+}{\eta},}
\end{array}
\end{equation}
where $h_+ = x_+ - x$ and $\tau = \hat{t} - t$, and notation~(\ref{SamarskiyNotation}) is used.

\medskip

Let us choose an invariant moving mesh of the form
\[
\frac{1}{\sqrt{\eta}}\left( \frac{\hat{x} - x}{\tau} - u \right) = 0
\]
or, equivalently,
\[
    \frac{\hat{x} - x}{\tau} = u.
\]
In the continuous limit it corresponds to the evolution of the spacial variable $x$ given as
\begin{equation} \label{dxlagr}
    \frac{dx}{dt} = u,
\end{equation}
where
\begin{equation} \label{dif}
\frac{d}{dt}= D_t + u D_x
\end{equation}
is the Lagrangian operator of total differentiation with respect to $t$.
The operator~(\ref{dif}) does not commute with $D_x$:
\[
\left[\frac{d}{dt}, D_x\right] \neq 0.
\]
Along with $\frac{d}{dt}$ we introduce two variables: a
``density'' $\rho$
\begin{equation} \label{var}
\rho = H + \eta,
\end{equation}
which in our case is just~$\eta$ because of the flat bottom;
and a new independent (mass) coordinate $s$ by means of {\it contact} transformation
\begin{equation} \label{s}
ds = \rho dx - \rho u dt,
\end{equation}
where $ds$ is a total differential form, i.e.,
\begin{equation} \label{tot}
\frac{\partial \rho}{\partial t} = - \frac{\partial(\rho u)}{\partial x}
\qquad
\text{or}
\qquad
\rho_t + (\rho u)_x =0.
\end{equation}
We also introduce the following operator of total
differentiation with respect to $s$:
\begin{equation} \label{D_s}
D_s = \frac{1}{\rho}D_x.
\end{equation}
The operators $\frac{d}{dt},D_s$ commute on the system
(\ref{Euler1}), (\ref{Euler2}):
\begin{equation} \label{oper}
\left[\frac{d}{dt}, D_s\right] = \left[D_t+ u D_x, \frac{1}{\rho}D_x
\right] = -\frac{1}{\rho^2}\left(\rho_t + (\rho u)_x\right)D_x =0.
\end{equation}
Then, in the new variables system~(\ref{Euler1}), (\ref{Euler2}) has the form
\begin{equation}\label{lagr_syst}
    \frac{d}{dt}\left(\frac {1}{\rho} \right) - D_s(u) =0,
    \qquad
    \frac{du}{dt} + \rho D_s(\rho-H) =0,
\end{equation}
where $H=0$ in the case under consideration.

Thus, we arrive at choosing the Lagrangian mass coordinates~$(t,s,u,\rho)$.
When discretizing equations, an important advantage of Lagrangian mass
coordinates is the possibility to construct invariant schemes on orthogonal
meshes~(see also \cite{dorodnitsyn2019shallow} for details).

Notice that a similar situation arises in the more general case of one-dimensional gas dynamics --- see~\cite{DORODNITSYN2019201}.

Invariant conservative schemes for the case of flat bottom
in Lagrangian coordinates and in Lagrangian mass coordinates were constructed by
the authors in the paper~\cite{dorodnitsyn2019shallow}.
Invariant schemes on adapting moving meshes for the flat bottom case were also considered in~\cite{bk:Bihlo_numeric[2012]}.
The schemes considered in~\cite{bk:Bihlo_numeric[2012]} do not possess the conservation law of energy.

\subsection{Parabolic bottom case}
\label{sec:parabolic}

In the parabolic bottom case, $H = \frac{1}{2}x^2$, on stencil~(\ref{14stencil})
one derives the following~10 invariants of the Lie algebra $X_1$, $X^3_2$, $X^3_3$, $X^3_4$ (see \#3 in Table~\ref{tab:HtableDelta})
\begin{equation}\label{quadInvs}
\begin{array}{c}
  {\displaystyle I^3_1 = \hat{t} - t,
  \quad
  I^3_2 = \frac{\hat{x}_+ - \hat{x}}{h_+},
  \quad
  I^3_3 = u_x,
  \quad
  I^3_4 = \hat{u}_x,}
  \\
  {\displaystyle
  I^3_5 = \frac{x^2 + 2 \eta}{h_+^2},
  \quad
  I^3_6 = \frac{x_+^2 + 2 \eta_+}{h_+^2},
  \quad
  I^3_7 = \frac{\hat{x}^2 + 2 \hat{\eta}}{h_+^2},
  \quad
  I^3_8 = \frac{\hat{x}_+^2 + 2 \hat{\eta}_+}{h_+^2},}
  \\
  {\displaystyle
  I^3_9 = \frac{(u + x)e^{2\tau} - 2 \hat{x}e^\tau - u + x}{(e^{2\tau} - 1) h_+},
  \quad
  I^3_{10} = \frac{2}{h_+}\left(
        \frac{\hat{x} - x}{\tau_1} - \frac{\hat{u} + u}{2}
        \right),}
\end{array}
\end{equation}
where it is denoted
\[
\tau_1
    = 2\,\frac{e^\tau - 1}{e^\tau + 1}
    = \frac{2 \sinh \tau}{1 + \cosh\tau} = \tau + O(\tau^3).
\]
Using invariants $I^3_1$ and $I^3_{10}$, one derives the following
Lagrangian-type moving mesh with flat time-layers
\begin{equation} \label{mesh_quad_lagr}
    \hat{t} - t = \tau = \textrm{const},
    \qquad
    \frac{\hat{x} - x}{\tau_1} = \frac{\hat{u} + u}{2}.
\end{equation}

In a similar way, one derives a Lagrangian-type moving mesh for $H = -\frac{1}{2}x^2$, namely
\begin{equation}
    {\displaystyle
    \tau = \textrm{const},
    \qquad
    \frac{\hat{x} - x}{\tau_2} = \frac{\hat{u} + u}{2},}
\end{equation}
where
\[
\tau_2
    = \frac{2 \sin \tau}{1 + \cos\tau} = \tau + O(\tau^3).
\]

Thus, we again arrive at Lagrangian-type coordinates.

\medskip

In the paper~\cite{dorodnitsyn2019shallow} the authors constructed
conservative invariant schemes for the shallow water equations
both in Lagrangian (potential) coordinates and in mass Lagrangian coordinates.
Here we slightly modify one of the mentioned schemes in order to construct
an invariant conservative scheme for the parabolic bottom.

Recall (see \cite{bk:KaptsovMeleshko_1D_classf[2018],dorodnitsyn2019shallow}) that in Lagrangian (potential) coordinates
the one-dimensional shallow water equations with the bottom $H(x)$ have the form
\begin{equation}\label{potential_eqns}
x_{tt} - 2\frac{x_{ss}}{x_s^3} - H^\prime(x) = 0,
\end{equation}
where $x=x(t,s)$ depends on the Lagrangian variable~$s$.

In the parabolic bottom case $H(x) = \pm\frac{x^2}{2}$, equation~(\ref{potential_eqns}) reads
\begin{equation}\label{potential_eqns_parabolic}
x_{tt} - 2\frac{x_{ss}}{x_s^3} \mp x = 0.
\end{equation}
Equation~(\ref{potential_eqns_parabolic}) possesses the conservation law of energy
\begin{equation} \label{CL1}
\frac{d}{dt}\left[ \frac{1}{x_s}+\frac{x_t^2}{2} \mp \frac{x^2}{2}\right]+  D_s\left[
\frac{x_t}{x_s^2}\right]=x_t\{x_{tt} - 2\frac{x_{ss}}{x_s^3} \mp x\}=0,
\end{equation}
and the conservation law of momentum
\begin{equation} \label{CL2}
\frac{d}{dt}\left[ x_t x_s\right]+  D_s\left[\frac{2}{x_s}-\frac{x_t^2}{2} \mp \frac{x^2}{2}\right]=x_s\{x_{tt}
- 2\frac{x_{ss}}{x_s^3} \mp x\}=0.
\end{equation}
The conservation law of mass is just
the symmetry of second derivatives relation, i.e.
\[
   \frac{d x_s}{dt} - D_s(x_t) = x_{ts} - x_{st} = 0.
\]

According to the classification provided in~\cite{bk:KaptsovMeleshko_1D_classf[2018]},
in addition to conservation laws of mass, energy and momentum,
there are two more conservation laws in the parabolic bottom case.

In case $H(x) = \frac{x^2}{2}$, the additional admitted generators and conservation laws are
\begin{enumerate}
\item
$e^t \frac{\partial}{\partial x}:$
\begin{equation} \label{CL1+quad}
\frac{d}{dt}\left[ e^t (x - x_t)\right]
    - D_s\left[e^t{x_s^{-2}}\right]
    = -e^t \{x_{tt} - 2\frac{x_{ss}}{x_s^3} - x\}=0;
\end{equation}

\item $e^{-t}\frac{\partial}{\partial x}:$
\begin{equation} \label{CL2+quad}
\frac{d}{dt}\left[ e^{-t} (x + x_t)\right]
    + D_s\left[e^{-t}{x_s^{-2}}\right]
    = e^{-t} \{x_{tt} - 2\frac{x_{ss}}{x_s^3} - x\}=0,
\end{equation}
\end{enumerate}
and their counterparts in mass coordinates are
\begin{enumerate}

\item
  $e^t\left(\frac{\partial}{\partial x} + \frac{\partial}{\partial u} - x \frac{\partial}{\partial \eta}\right):$
  \begin{equation}\label{CL1+quad_mass}
    \frac{d}{dt}(e^t (x - u)) - \frac{1}{2}D_s(e^t \rho^2)
    = e^t \left\{\frac{dx}{dt} - u\right\}
        -e^t\left\{\frac{du}{dt} + \rho \left(\rho - \frac{x^2}{2}\right)_s\right\}  = 0;
  \end{equation}

\item
  $e^{-t}\left(\frac{\partial}{\partial x} - \frac{\partial}{\partial u} - x \frac{\partial}{\partial \eta}\right):$
  \begin{equation}\label{CL2+quad_mass}
    \frac{d}{dt}(e^{-t} (x + u)) + \frac{1}{2}D_s(e^{-t}\rho^2)
    = e^{-t} \left\{\frac{dx}{dt} - u\right\}
        +e^{-t}\left\{\frac{du}{dt} + \rho \left(\rho - \frac{x^2}{2}\right)_s\right\}  = 0,
  \end{equation}
  where $\rho = \eta + H = \eta + \frac{x^2}{2}$.
\end{enumerate}

In case $H(x) = -\frac{x^2}{2}$, the conservation laws are
\begin{enumerate}
\item
$\sin t \frac{\partial}{\partial x}:$
\begin{equation} \label{CL1-quad}
\frac{d}{dt}\left[ x\cos t - x_t \sin t \right]
    - D_s\left[{x_s^{-2}}\sin t\right]
    = -\sin t \{x_{tt} - 2\frac{x_{ss}}{x_s^3} - x\}=0;
\end{equation}

\item
$\cos t \frac{\partial}{\partial x}:$
\begin{equation} \label{CL2-quad}
\frac{d}{dt}\left[ x\sin t + x_t \cos t \right]
    + D_s\left[{x_s^{-2}}\cos t\right]
    = \cos t \{x_{tt} - 2\frac{x_{ss}}{x_s^3} - x\}=0;
\end{equation}
\end{enumerate}
and their counterparts in mass coordinates are
\begin{enumerate}

\item
  $\sin t \frac{\partial}{\partial x} + \cos t \frac{\partial}{\partial u}
  + x \sin t \frac{\partial}{\partial \eta}:$
  \begin{equation}\label{CL1-quad_mass}
    \frac{d}{dt}\left[ x\cos t - u \sin t \right]
    - D_s\left[\rho^2 \sin t\right] = 0;
  \end{equation}

\item
  $\cos t \frac{\partial}{\partial x} - \sin t \frac{\partial}{\partial u}
  + x \cos t \frac{\partial}{\partial \eta}:$
  \begin{equation}\label{CL2-quad_mass}
    \frac{d}{dt}\left[ x\sin t + u \cos t \right]
        + D_s\left[\rho^2 \cos t\right] = 0,
  \end{equation}
  where $\rho = \eta + H = \eta - \frac{x^2}{2}$.
\end{enumerate}

Let us begin with a scheme for the case $H(x) = \frac{x^2}{2}$.
To obtain a conservative invariant scheme we extend the scheme constructed in~\cite{dorodnitsyn2019shallow}
\begin{equation} \label{schemeLagrMain}
  \displaystyle
  \def\arraystretch{1.75}
  \begin{array}{c}
        \displaystyle
        x_{t\check{t}}
      + \frac{1}{h^s_-} \left(
        (\hat{x}_s \check{x}_s)^{-1}
        - (\hat{x}_{\bar{s}} \check{x}_{\bar{s}})^{-1}
      \right)
      = 0,
      \\
      \displaystyle
      \tau_+ = \tau_-, \qquad
      h^s_+ = h^s_-
  \end{array}
\end{equation}
as follows,
\begin{equation}
  \displaystyle
  \def\arraystretch{1.75}
  \begin{array}{c}
        \displaystyle
        F = x_{t\check{t}}
      + \frac{1}{h^s_-} \left(
        (\hat{x}_s \check{x}_s)^{-1}
        - (\hat{x}_{\bar{s}} \check{x}_{\bar{s}})^{-1}
      \right) - x \phi(\tau)
      = 0,
      \\
      \displaystyle
      \Omega: \qquad
      \tau_+ = \tau_-, \qquad
      h^s_+ = h^s_-,
  \end{array}
\end{equation}
where $\phi(\tau)$ is an unknown function which tends to $1$ as $\tau$ tends to zero,
and
\[
h^s_+ = s_+ - s = h^s,
\qquad
h^s_- = s - s_- = h^s.
\]

Applying the difference
variational Euler operator~(see~\cite{bk:Dorodnitsyn[2011],bk:Dorod_Hamilt[2011]}
and Section~\ref{sec:scheme_Euler_arb} for details)
\begin{equation} \label{EulerOp}
\mathcal{E}_x =
    \sum_{k=-\infty}^{+\infty}
    \sum_{l=-\infty}^{+\infty}
    {\dtp{S}}^k{\dsp{S}}^l
    \frac{\partial}{\partial x^{n-k}_{m-l}}
\end{equation}
on the uniform mesh $\Omega$
\[
    \mathcal{E}_x(e^{\pm t} F)|_{\Omega}
    = \frac{1}{\tau^2} e^{\pm t} (\tau^2 \phi - e^\tau - e^{-\tau} + 2)
    = 0,
\]
one obtains the solution, namely
\[
    \phi(\tau) = \frac{2 (\cosh \tau - 1)}{\tau^2} = 1 + \frac{\tau^2}{12} + O(\tau^4).
\]
Thus, the invariant scheme becomes
\begin{equation} \label{scheme_sq_lagr}
  \displaystyle
  \def\arraystretch{1.75}
  \begin{array}{c}
        \displaystyle
        x_{t\check{t}}
      + \dsm{D}\left(
        \frac{1}{\hat{x}_s \check{x}_s}
      \right) - \frac{2 (\cosh \tau - 1)}{\tau^2} x
      = 0,
      \\
      \displaystyle
      \qquad
      \tau_+ = \tau_-, \qquad
      h^s_+ = h^s_-.
  \end{array}
\end{equation}

\begin{remark}
\label{rem:parabHequiv}
It is essential for numerical computations that for~$H(x) = \frac{\beta}{2} (x - c)^2$, where $\beta = \text{const} >0$ and $c = \text{const}$,
the conservation law multipliers are~$e^{\pm\sqrt{\beta}t}$ and
\[
    \phi(\tau) = \frac{2 (\cosh(\sqrt{\beta} \tau) - 1)}{(\sqrt{\beta} \tau)^2}.
\]
Similar changes should be done for~$H(x) = -\frac{\beta}{2}(x - c)^2$.
\end{remark}

Scheme~(\ref{scheme_sq_lagr}) possesses the conservation laws of mass and energy
\begin{equation}\label{csmassLagrMass}
  \dtm{D}(\hat{x}_s) - \dsm{D}(x_t^+) = 0,
\end{equation}
\begin{equation} \label{CLparabEnergy}
    \dtm{D}(
        x_t^2 + x_s^{-1} + \hat{x}_s^{-1}
    )
    + \dsm{D}\left(
        (x_t^+ + \check{x}_t^+)
        (\hat{x}_s \check{x}_s)^{-1}
        -  \frac{2 (\cosh \tau - 1)}{\tau^2}  x\hat{x}
    \right) = 0,
  \end{equation}
and the difference analogues of conservation laws (\ref{CL1+quad}) and (\ref{CL2+quad}), i.e.,
\[
{\dtm{D}}\left(
    x \frac{e^{\hat{t}} - e^t}{\tau}
    - e^t {x}_t
\right)
    -\dsm{D}\left(
        e^t \frac{1}{\hat{x}_s \check{x}_s}
    \right)
    =
    -e^t \left\{
        x_{t\check{t}}
      + \dsm{D}\left(
        \frac{1}{\hat{x}_s \check{x}_s}
      \right) - \frac{2 (\cosh \tau - 1)}{\tau^2} x
    \right\} = 0,
\]
\[
{\dtm{D}}\left(
    x \frac{e^{-t} - e^{-\hat{t}}}{\tau}
    - e^{-t} {x}_t
\right)
    +\dsm{D}\left(
        e^{-t} \frac{1}{\hat{x}_s \check{x}_s}
    \right)
    =
    e^{-t} \left\{
        x_{t\check{t}}
      + \dsm{D}\left(
        \frac{1}{\hat{x}_s \check{x}_s}
      \right) - \frac{2 (\cosh \tau - 1)}{\tau^2} x
    \right\} = 0.
\]
Notice that $\frac{1}{\tau}(e^{\hat{t}} - e^{t}) = e^t + O(\tau)$ as in continuous case $(e^t)^\prime = e^t$.

\medskip

By a similar procedure one arrives at the following scheme for
the case $H(x) = -\frac{x^2}{2}$
\begin{equation} \label{scheme_sq_lagr2}
  \displaystyle
  \def\arraystretch{1.75}
  \begin{array}{c}
        \displaystyle
        x_{t\check{t}}
      + \dsm{D}\left(
        \frac{1}{\hat{x}_s \check{x}_s}
      \right) - \frac{2 (\cos \tau - 1)}{\tau^2} x
      = 0,
      \\
      \displaystyle
      \qquad
      \tau_+ = \tau_-, \qquad
      h^s_+ = h^s_-,
  \end{array}
\end{equation}
where
\[
\frac{2 (\cos \tau - 1)}{\tau^2} = -1 + \frac{\tau^2}{12} + O(\tau^4).
\]
The scheme possesses the conservation laws of mass~(\ref{csmassLagrMass}), energy
\[
    \dtm{D}(
        x_t^2 + x_s^{-1} + \hat{x}_s^{-1}
    )
    + \dsm{D}\left(
        (x_t^+ + \check{x}_t^+)
        (\hat{x}_s \check{x}_s)^{-1}
        -  \frac{2 (\cos \tau - 1)}{\tau^2}  x\hat{x}
    \right) = 0,
\]
and two additional conservation laws:
\[
{\dtm{D}}\left(
    {x}_t \sin t
    - x \frac{\sin \hat{t} - \sin t}{\tau}
\right)
    +\dsm{D}\left(
        \sin t \frac{1}{\hat{x}_s \check{x}_s}
    \right)
    =
    \sin t \left\{
        x_{t\check{t}}
      + \dsm{D}\left(
        \frac{1}{\hat{x}_s \check{x}_s}
      \right) - \frac{2 (\cos \tau - 1)}{\tau^2} x
    \right\} = 0,
\]
\[
{\dtm{D}}\left(
    {x}_t \cos t
    - x \frac{\cos \hat{t} - \cos t}{\tau}
\right)
    +\dsm{D}\left(
        \cos t \frac{1}{\hat{x}_s \check{x}_s}
    \right)
    =
    \cos t \left\{
        x_{t\check{t}}
      + \dsm{D}\left(
        \frac{1}{\hat{x}_s \check{x}_s}
      \right) - \frac{2 (\cos \tau - 1)}{\tau^2} x
    \right\} = 0.
\]

\medskip

In order to transform the three-layer scheme into a two-layer one
we apply the same technique that was used in~\cite{dorodnitsyn2019shallow}.
It consists of a special approximation of the ``state equation''~$p = \rho^2$.
Then, scheme~(\ref{scheme_sq_lagr}) can be represented
in hydrodynamic variables on two time layers as follows
\begin{equation} \label{mass_lagr_sq}
    \def\arraystretch{2}
    \begin{array}{c}
    \displaystyle
    {\dtm{D}}\left( \frac{1}{\rho} \right) - \dsm{D}\left(
        \frac{u^+ + \check{u}^+}{2}
    \right) = 0,
    \\
    \displaystyle
      {\dtm{D}}(u) + \dsm{D}\left(Q\right) - \frac{2 (\cosh \tau - 1)}{\tau^2} x = 0,
    \\
    \displaystyle
    \check{x}_s + x_s = \frac{1}{\sqrt{\check{p}}} + \frac{1}{\sqrt{p}} = \frac{2}{\check{\rho}},
    \\
    x_t = u, \quad
    h^s_+ = h^s_-, \quad
    \tau_+ = \tau_-,
    \end{array}
\end{equation}
where $Q$ is given by the relation
\begin{equation}\label{flux_Q}
  \frac{1}{Q} = \frac{4}{\rho\check{\rho}}
    - \frac{2}{\sqrt{p}}\left( \frac{1}{\rho} + \frac{1}{\check{\rho}} \right) + \frac{1}{p},
\end{equation}
and the equation
\[
    \check{x}_s + x_s = \frac{1}{\sqrt{\check{p}}} + \frac{1}{\sqrt{p}} = \frac{2}{\check{\rho}}
\]
approximates the equation~$p = \rho^2$.

\begin{figure}[h]
\centering
\includegraphics[scale=0.75]{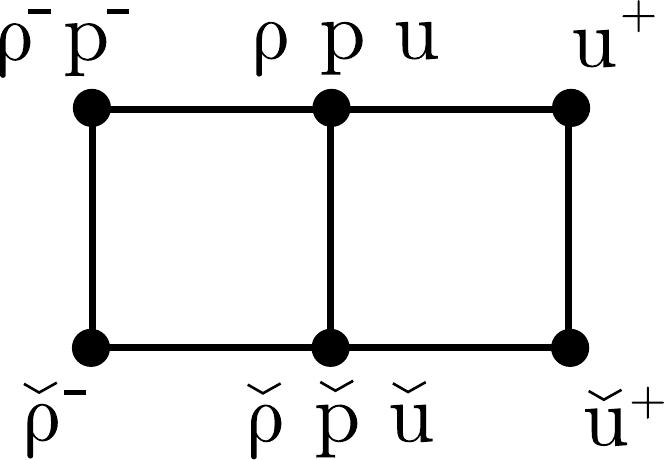}
\caption{Two-layer scheme difference template}
\label{fig:template_t2}
\end{figure}

Thus, a decrease in the number of time layers of the scheme is achieved by
increasing the number of equations in the system.
The two time-layer template for scheme~(\ref{mass_lagr_sq}) is depicted in Figure~\ref{fig:template_t2}.

Notice that $x = \frac{1}{2}\dsm{D}(xx_+)$ on the uniform mesh,
so the second equation of~(\ref{mass_lagr_sq}) can be rewritten in the following divergent form
\[
{\dtm{D}}(u) + \dsm{D}\left(Q - \frac{\cosh \tau - 1}{\tau^2} x x_+ \right) = 0.
\]

Similarly, scheme~(\ref{scheme_sq_lagr2}) can be rewritten as follows
\begin{equation} \label{mass_lagr_sq2}
    \def\arraystretch{2}
    \begin{array}{c}
    \displaystyle
    {\dtm{D}}\left( \frac{1}{\rho} \right) - \dsm{D}\left(
        \frac{u^+ + \check{u}^+}{2}
    \right) = 0,
    \\
    \displaystyle
      {\dtm{D}}(u) + \dsm{D}\left(Q\right) - \frac{2 (\cos \tau - 1)}{\tau^2} x = 0,
    \\
    \displaystyle
    \check{x}_s + x_s = \frac{1}{\sqrt{\check{p}}} + \frac{1}{\sqrt{p}} = \frac{2}{\check{\rho}},
    \\
    x_t = u, \quad
    h^s_+ = h^s_-, \quad
    \tau_+ = \tau_-,
    \end{array}
\end{equation}
where $Q$ is given by~(\ref{flux_Q}).

The difference analogues of conservation laws~(\ref{CL1+quad_mass}) and~(\ref{CL2+quad_mass})
of scheme~(\ref{mass_lagr_sq}) are as follows,
\begin{equation}
{\dtm{D}}\left(
    x \frac{e^{\hat{t}} - e^t}{\tau}
    - e^t u
\right)
    -\dsm{D}\left(
        e^t Q
    \right)
    = 0,
\end{equation}
\begin{equation}
{\dtm{D}}\left(
    x \frac{e^{-t} - e^{-\hat{t}}}{\tau}
    - e^{-t} u
\right)
    +\dsm{D}\left(
        e^{-t} Q
    \right)
    = 0.
\end{equation}

The difference analogues of conservation laws~(\ref{CL1-quad_mass}) and~(\ref{CL2-quad_mass})
of scheme~(\ref{mass_lagr_sq2}) are
\begin{equation}
{\dtm{D}}\left(
    x \frac{\sin \hat{t} - \sin t}{\tau}
    - u \sin t
\right)
    -\dsm{D}\left(
        Q \sin t
    \right) = 0,
\end{equation}
\begin{equation}
{\dtm{D}}\left(
    x \frac{\cos \hat{t} - \cos t}{\tau}
    - u \cos t
\right)
    -\dsm{D}\left(
        Q \cos t
    \right)
    = 0.
\end{equation}

\begin{remark}
In the present section we have extended scheme~(\ref{schemeLagrMain})
to the case of a parabolic bottom.
Scheme~(\ref{schemeLagrMain}) can be successfully extended to the case of
the linear bottom
\[
H(x)=C_1 x + C_2
\]
as well:
\[
  \displaystyle
  \def\arraystretch{1.75}
  \begin{array}{c}
    \displaystyle
      x_{t\check{t}}
        + \frac{1}{h^s_-} \left(
        (\hat{x}_s \check{x}_s)^{-1}
        - (\hat{x}_{\bar{s}} \check{x}_{\bar{s}})^{-1}
      \right)
      - C_1
      = 0,
      \\
      \displaystyle
      \tau_+ = \tau_-,
      \qquad
      h^s_+ = h^s_-
  \end{array}
\]
It was stated in~\cite{dorodnitsyn2019shallow} that the latter scheme
is related to~(\ref{schemeLagrMain})
by the following transformation
\begin{equation}\label{linearbtmTrLagr}
x = \tilde{x} + \frac{C_1}{2} t\hat{t},
\qquad
t = \tilde{t},
\qquad
s = \tilde{s}.
\end{equation}
\end{remark}

\begin{remark}
\label{rem:H-4/3}
According to the content of Table~\ref{tab:Htable}
and Table~\ref{tab:HtableDelta},
no additional (differential or finite-difference) conservation laws
occur in Eulerian coordinates in the case of the bottom
\[
H(x) = x^c,
\qquad
c = \textrm{const} \neq 0, 1, 2.
\]
In contrast to Eulerian coordinates, in Lagrangian coordinates an additional
differential conservation law occurs for~$c={-4/3}$~(see~\cite{bk:KaptsovMeleshko_1D_classf[2018]} for details).
In Lagrangian coordinates, a finite-difference analogue of that conservation law may
exist for some difference scheme.
We leave this very particular case out of consideration here.
\end{remark}

\section{Conservative invariant scheme for an arbitrary bottom shape in Eulerian coordinates}
\label{sec:scheme_Euler_arb}

An analysis of finite-difference invariants~(\ref{flatInvs})
and~(\ref{quadInvs}) for the flat and parabolic bottom topography
indicates that it is unlikely possible to construct
invariant schemes on orthogonal meshes in Eulerian coordinates.
In contrast, Lagrangian coordinates allow one to construct
invariant schemes on \emph{uniform orthogonal} meshes.
In~\cite{dorodnitsyn2019shallow} the authors constructed
invariant conservative schemes in Lagrangian coordinates for the
flat bottom topography. However, the authors faced some difficulties in
constructing schemes for an arbitrary bottom topography. Apparently, in
that case it is only possible to construct schemes that preserve
mass and either energy or momentum.

Here we show that in Eulerian coordinates it is possible to
construct schemes that possess conservation laws of mass, momentum
and energy simultaneously.

Notice that according to the content of Table~\ref{tab:Htable}, the
only generator admitted by the equations is the time shift
generator~$\frac{\partial}{\partial t}$ and, therefore, the only
restriction imposed by the invariance is that the scheme can not
explicitly depend on the variable~$t$. Thus, we can only stay
focused on constructing conservative schemes. In order to construct
such schemes, we use the finite-difference analogue of the
direct method~\cite{[Bluman1],bk:BlumanCheviakovAnco}. The
difference analogue of the direct method is useful both in obtaining
conservation laws of known difference schemes, and in constructing
new conservative difference schemes provided with certain preliminary
assumptions on their forms. The approach can be fairly effective
when considering polynomial schemes or schemes that are some
rational functions defined on the chosen difference stencil. The
polynomial schemes approach was successfully applied by the authors
in~\cite{dorodnitsyn2019shallow,bk:ChevDorKap2020}. The method is
also used in Section~\ref{sec:parabolic} above for
schemes~(\ref{scheme_sq_lagr}) and~(\ref{scheme_sq_lagr2}).

The key idea of the direct method is to apply the variational Euler
operator~(\ref{EulerOp}) on a uniform orthogonal mesh to some
difference approximation (a family of schemes) with undetermined
coefficients which should be expressed in a divergent form.
Under the action of operator~(\ref{EulerOp}) any finite-difference
divergent expression vanishes that allows one to find the
coefficients and since  to obtain specific schemes and their
conservation laws.

\medskip

Consider the shallow water equations~(\ref{Euler1}), (\ref{Euler2}) and their
conservation law of energy~(\ref{Eul1})
\[
D_t\left[u^2(\eta +H)+ \eta^2\right]
+ D_x\left[u(\eta +H(x))(u^2 +2\eta)\right] = 0
\]
in the following form
\begin{equation} \label{MuEqn}
\mu_1 (\eta_t +((\eta+H(x))u)_x) - \mu_2 (u_t+uu_x+\eta_x) = 0,
\end{equation}
where
\begin{equation}
\mu_1 = \eta +\frac{u^2}{2}
\qquad
\text{and}
\qquad
\mu_2 = (\eta + H(x)) u
\end{equation}
are the conservation law multipliers (or integrating multipliers).

We approximate equations~(\ref{Euler1}), (\ref{Euler2}) on the 4-point stencil~(\ref{14stencil})
by finite-difference polynomials
\begin{equation} \label{constr:F1}
  F_1[u,\eta,H] = \dtp{D}(p_1 \eta^n_m + (1 - p_1) \eta^n_{m+1})
  + \dhp{D}\left(
        \sum_{0 \leqslant k,p \leqslant 1} w_{k+1,p+1} (\eta^{n+k}_m + H(x_m)) u^{n+p}_m
  \right) = 0,
\end{equation}
and
\begin{multline} \label{constr:F2}
  F_2[u,\eta] = \dtp{D}(q_1 u^n_m + (1 - q_1) u^{n}_{m+1})
  \\
  + \dhp{D}\left(
        \frac{1}{2}(z_{11} (u^n_m)^2 + z_{12} u^n_m u^{n+1}_m + z_{22} (u^{n+1}_m)^2)
        + q_2 \eta^n_m + (1 - q_2) \eta^{n+1}_m
  \right) = 0,
\end{multline}
and the integrating multipliers~$\mu_1$ and $\mu_2$ by the following expressions
\begin{equation} \label{M1}
  M_1[u,\eta] = \sum_{0 \leqslant k,l \leqslant 1} B_{k+1,l+1} \eta^{n+k}_{m+l}
    + \frac{1}{2} \sum_{0 \leqslant k,l,p,q \leqslant 1} a_{k+1,l+1,p+1,q+1} u^{n+k}_{m+l} u^{n+p}_{m+q},
\end{equation}
\begin{equation} \label{M2}
    M_2[u,\eta,H] = \sum_{0 \leqslant k,l,p,q \leqslant 1} b_{k+1,l+1,p+1,q+1} (\eta^{n+k}_{m+l} + H(x_m)) u^{n+p}_{m+q},
\end{equation}
where $p_i$, $q_i$, $B_{ij}$, $z_{ij}$, $w_{ij}$, $a_{ijkl}$ and $b_{ijkl}$
are some constant undetermined coefficients.

It follows from the form of the chosen approximation that the constants are related by
\begin{equation} \label{constrels}
    \displaystyle
    \sum_{i,j,k,l} a_{ijkl} = 1,
    \qquad
    \sum_{i,j,k,l} b_{ijkl} = 1,
    \qquad
    \sum_{i,j} B_{ij} = 1,
    \qquad
    \sum_{i,j} z_{ij} = 1.
\end{equation}
For example the last relation means
\[
z_{11} (u^n_m)^2 + z_{12} u^n_m u^{n+1}_m + z_{22} (u^{n+1}_m)^2 \sim u^2.
\]

Obviously, approximation~(\ref{constr:F1}),~(\ref{constr:F2})
for any set of coefficients admits the
generator~$\frac{\partial}{\partial t}$, i.e.,
the scheme is a invariant one.

According to the direct method, we require
\begin{equation}\label{EL}
  \mathcal{E}_u (M_1 F_1 + M_2 F_2) \equiv 0,
  \qquad
  \mathcal{E}_\eta (M_1 F_1 + M_2 F_2) \equiv 0,
\end{equation}
where
\begin{equation}\label{DirectCLenergy}
M_1 F_1 + M_2 F_2 = 0
\end{equation}
approximates the conservation law of energy~(\ref{MuEqn}).
Notice that equations~(\ref{constr:F1}),~(\ref{constr:F2}) are conservative by construction,
i.e., the  conservation laws of momentum and mass hold.

Considering~(\ref{EL}) and taking~(\ref{constrels}) into account,
we obtain several sets of relations on the coefficients~$p_i, q_i, B_{ij}, z_{ij}, w_{ij}, a_{ijkl}, b_{ijkl}$.
It can be shown by means of~(\ref{constrels}) that these sets are not independent.
Thus, without loss of generality, we consider only one set of relations, namely:
\begin{equation} \label{indetCoeffs}
\def\arraystretch{1.5}
\begin{array}{c}
B_{12} = B_{22},
\quad
a_{1112} = -a_{1211},
\quad
a_{1121} = -a_{2111},
\quad
a_{1122} = -a_{2211},
\\
a_{1212} = 2 b_{2111},
\quad
a_{1221} = -a_{2112},
\quad
a_{1222} = 2 B_{22} z_{12} - a_{2212},
\quad
a_{2122} = -a_{2221},
\\
a_{2222} = 4 B_{22} w_{11} - 2 B_{22} z_{12},
\quad
b_{1111} = 2 B_{22} w_{11},
\quad
b_{1121} = 2 B_{22} w_{11} - B_{22} z_{12},
\\
b_{2121} = B_{22} z_{12} + b_{2111},
\quad
q_{1} = 1,
\quad
q_{2} = \frac{1}{2},
\quad
w_{12} = w_{11}- \frac{1}{2}z_{12},
\quad
w_{21} = \frac{b_{2111}}{2 B_{22}},
\\
w_{22} = \frac{B_{22} z_{12} + b_{2111}}{2 B_{22}},
\quad
z_{11} = \frac{b_{2111}}{B_{22}},
\quad
z_{22} = 2 w_{11} - z_{12},
\end{array}
\end{equation}
where $B_{22}$, $a_{1211}$, $a_{2111}$, $a_{2112}$, $a_{2211}$, $a_{2212}$, $a_{2221}$, $b_{2111}$, $w_{11}$ and $z_{12}$ are arbitrary provided relations~(\ref{constrels}),
and the rest coefficients are zero.

Substituting~(\ref{indetCoeffs}) into~(\ref{constr:F1})--(\ref{M2})
and then expanding the resulting expressions into series, we derive that
\begin{equation} \label{approxSeries}
\begin{array}{c}

F_1 \sim \eta_t + \frac{2 B_{22} w_{11} + b_{2111}}{B_{22}} \left(
        (\eta + H) u
    \right)_x,
\qquad
F_2 \sim u_t + \eta_x + \frac{2 B_{22} w_{11} + b_{2111}}{B_{22}} u u_x,
\\
M_1 \sim (2 B_{22} w_{11} + b_{2111}) u^2 + 2 B_{22} \eta,
\qquad
M_2 \sim 2 u (\eta + H) (2 B_{22} w_{11} + b_{2111}).
\end{array}
\end{equation}
Analyzing~(\ref{approxSeries}), we conclude that it should be
\begin{equation}
w_{11} + b_{2111} = B_{22} = \frac{1}{2}.
\end{equation}
Finally, the coefficients $z_{12}$ and $w_{11}$ remain arbitrary,
and the multipliers $M_1$ and $M_2$ are
\begin{equation}
M_1 = \frac{1}{2}\left[
        (2 w_{11} - z_{12}) (u^{n+1}_{m+1})^2
        + z_{12} u^n_{m+1} u^{n+1}_{m+1}
        + (1 - 2 w_{11}) (u^n_{m+1})^2
        + \eta^n_{m+1}
        + \eta^{n+1}_{m+1}
    \right],
\end{equation}
\begin{equation}
\begin{array}{c}
M_2 = \frac{1}{2}\Big[
        (u^n_m + u^{n+1}_m) H
        + \left((
            1 + z_{12} - 2 w_{11}) \eta^{n+1}_m
            + (2  w_{11} - z_{12}) \eta^n_m
        \right) u^{n+1}_m
        \\
        + \left(
            (1 - 2 w_{11}) \eta^{n+1}_m + 2 w_{11} \eta^n_m
        \right) u^n_m
    \Big].
\end{array}
\end{equation}
The conservation law of energy~(\ref{DirectCLenergy}) possesses the form
\begin{multline} \label{scheme_longform}
  \frac{1}{2}\Bigg\{
    \dtp{D}\left(
        {\eta_+}^2 + (\eta + H) u^2 + w_{11} \tau \dhp{D}(H u^3)
    \right)
    +
    \dhp{D}\Bigg(
        (\hat{\eta} + H)
        \frac{u + \hat{u}}{2}
        (u^2 + \eta + \hat{\eta})
        + h \dtp{D}(\eta) u^2
        \\
        - 2\tau^2\!\left[
            (u + \hat{u})^2 w_{11}^2
            - (u + \hat{u}) \hat{u} w_{11}z_{12}
            + \frac{1}{4} \hat{u}^2 z_{12}^2
        \right]\dtp{D}(\eta)\!\dtp{D}(u)
        + \Theta_1 w_{11}
        + \Theta_2 z_{12}
    \Bigg)
  \Bigg\} = 0,
\end{multline}
where
\[
\Theta_1 = \tau u \hat{u} \dtp{D}(u) H
    - \left[
        \hat{\eta}^2
        + (2u^2 - \hat{u}^2)\hat{\eta} - (u^2 + \eta)\eta
    \right] (u + \hat{u})
    + h\tau(u + \hat{u})\dtp{D}(\eta)\!\dtp{D}(u),
\]
\[
\displaystyle
\Theta_2 = \frac{\hat{u}}{2}\left[
        (u^2 - \hat{u}^2) H
        + \hat{\eta}^2
        + (2u^2 - \hat{u}^2)\hat{\eta} - (u^2 + \eta)\eta
    \right]
    - h\tau\hat{u}\dtp{D}(\eta)\!\dtp{D}(u).
\]

Thus, the two-parametric ($z_{12}$ and $w_{11}$) family of
conservative invariant schemes has been obtained.

\medskip

Here we consider two particular cases of scheme~(\ref{constr:F1}), (\ref{constr:F2}).
\begin{enumerate}

\item
To eliminate the terms depending on $\tau$, in~(\ref{scheme_longform}) we set $w_{11} = 0$ and $z_{12} = 0$.
As a result, we derive the following scheme
\begin{equation} \label{scm_arb}
    \def\arraystretch{1.75}
    \large
    \begin{array}{c}
        \dtp{D}(\eta_+) + \frac{1}{2}\dhp{D}\left( (\hat{u} + u) (\hat{\eta} + H) \right) = 0,
        \\
        \dtp{D}(u) + \frac{1}{2}\dhp{D}\left( u^2 + \hat{\eta} + \eta\right) = 0.
    \end{array}
\end{equation}
The latter scheme possesses the conservation law of energy
\begin{multline*} 
  \frac{1}{2}(u_+^2 + \eta_+ + \hat{\eta}_+) \left\{
        \dtp{D}(\eta_+) + \dhp{D}\left( \frac{u + \hat{u}}{2} (\hat{\eta} + H) \right)
  \right\}
  \\
  +
  \frac{(\hat{\eta} + H)(u + \hat{u})}{2} \left\{
       \dtp{D}(u) + \frac{1}{2}\dhp{D}\left( u^2 + \eta + \hat{\eta}\right)
  \right\}
  \\
  =
  \frac{1}{2}\left\{
    \dtp{D}\left(
        u^2 (\eta + H) + \eta_+^2
    \right)
    +
    \dhp{D}\left(
        (\hat{\eta} + H)
        \frac{u + \hat{u}}{2}
        (u^2 + \eta + \hat{\eta})
        + h \dtp{D}(\eta) u^2
    \right)
  \right\} = 0
\end{multline*}
as well as conservation laws of mass and momentum.
Scheme~(\ref{scm_arb}) is the simplest scheme of the considered family
as it corresponds to zero-valued parameters~$w_{11}$ and~$z_{12}$.

\item
To give the difference system a slightly more symmetric form
we can set~$w_{11}=\frac{1}{2}$ and~$z_{12} = 1$
and derive the scheme
\begin{equation} \label{scm_sym}
    \def\arraystretch{1.75}
    \large
    \begin{array}{c}
        \dtp{D}(\eta_+)
        + \frac{1}{2}\dhp{D}\left(
            \eta u + \hat{\eta}\hat{u} + (\hat{u} + u) H
        \right) = 0,
        \\
        \dtp{D}(u) + \frac{1}{2}\dhp{D}\left( u\hat{u} + \hat{\eta} + \eta\right) = 0.
    \end{array}
\end{equation}
This scheme possesses the conservation laws of mass and momentum
and the conservation law of energy
\begin{multline}\label{scm_sym_arb_cl}
  \frac{1}{2}(\hat{u}_+ u_+ + \hat{\eta}_+ + \eta_+) \left\{
        \dtp{D}(\eta_+)
        + \frac{1}{2}\dhp{D}\left(
            \eta u + \hat{\eta}\hat{u} + (\hat{u} + u) H
        \right)
  \right\}
  \\
  +
  \frac{1}{2}
  \left(\hat{u}\hat{\eta} + u\eta + (\hat{u} + u) H\right)
  \left\{
       \dtp{D}(u) + \frac{1}{2}\dhp{D}\left( u\hat{u} + \hat{\eta} + \eta\right)
  \right\}
  \\
  =
  \frac{1}{2}
    \dtp{D}\left(
        u^2 (\eta + H) + \eta_+^2
    \right)
    + \frac{1}{4}\dhp{D}\Big\{
    (u\hat{u} + \hat{\eta} + \eta)
    \left(
        \hat{u}\hat{\eta} + u\eta
        + (\hat{u} + u) H
    \right)
        + 2hu\hat{u}\!\dtp{D}(\eta)
    \Big\}
  = 0.
\end{multline}
\end{enumerate}
Notice that the schemes obtained are not unique.

\begin{remark}
In terms of variable $\rho = \eta + H$ one can rewrite equations~(\ref{Euler1}), (\ref{Euler2}) as follows
\[
\eta_t + (\rho u)_x = 0,
\]
\[
u_t +uu_x + \rho_x = H^\prime.
\]
Then, scheme~(\ref{scm_arb}) becomes
\[
    \def\arraystretch{1.75}
    \large
    \begin{array}{c}
        \dtp{D}(\rho_+) + \frac{1}{2}\dhp{D}\left( (\hat{u} + u) \hat{\rho} \right) = 0,
        \\
        \dtp{D}(u) + \frac{1}{2}\dhp{D}\left( u^2 + \hat{\rho} + \rho\right) = \dhp{D}(H),
    \end{array}
\]
and scheme~(\ref{scm_sym}) possesses the form
\[
    \def\arraystretch{1.75}
    \large
    \begin{array}{c}
        \dtp{D}(\rho_+)
        + \frac{1}{2}\dhp{D}\left(
            \rho u + \hat{\rho}\hat{u}
        \right) = 0,
        \\
        \dtp{D}(u) + \frac{1}{2}\dhp{D}\left( u\hat{u} + \hat{\rho} + \rho\right) = \dhp{D}(H).
    \end{array}
\]
The latter forms of the equations look simpler than the original ones.
\end{remark}

\begin{remark}
In case of a flat bottom, the shallow water equations admit the generator~$\frac{\partial}{\partial x}$~(which is a particular case of the generator~$X^1_\infty$ --- see Table~\ref{tab:Htable}), and the corresponding conservation law (of momentum) is
\begin{equation} \label{CLnone}
  u (\eta_t + \eta u_x + u \eta_x) + \eta (u_t + u u_x + \eta_x)
  = D_t(u \eta) + D_x\left( \eta u^2 + \frac{\eta^2}{2} \right) = 0.
\end{equation}
It can be shown by a procedure similar to that described above
that there are no polynomial schemes of form~(\ref{constr:F1}), (\ref{constr:F2})
possessing the difference analogue of conservation law~(\ref{CLnone}) with
integrating multipliers linear by~$u$ and~$\eta$.
\end{remark}


\section{Numerical implementation of invariant
conservative schemes in Eulerian coordinates}

In the present section we consider scheme~(\ref{scm_sym})
and analyze some of its numerical properties.
Then we compare the scheme with a modified scheme which
does not preserve energy.

\subsection{Implementation of scheme~(\ref{scm_sym})}
\label{sec:Euler_scm_impl}

Let us rewrite scheme~(\ref{scm_sym}) in index form as follows
\begin{equation} \label{scm_sym_index}
    \def\arraystretch{1.75}
    \begin{array}{c}
        \eta^{n+1}_{m+1} - \eta^{n}_{m+1} + {a}\big[
            \eta^{n}_{m+1} u^{n}_{m+1} + \eta^{n+1}_{m+1} u^{n+1}_{m+1}
            - \eta^{n}_{m} u^{n}_{m} - \eta^{n+1}_{m} u^{n+1}_{m}
            \\
            + (u^{n+1}_{m+1} + u^n_{m+1}) H(x_{m+1})
            - (u^{n+1}_{m} + u^n_{m}) H(x_{m})
        \big] = 0,
        \\
        u^{n+1}_m - u^n_m + {a}\left[
            u^n_{m+1}u^{n+1}_{m+1} - u^n_m u^{n+1}_m
            + \eta^{n+1}_{m+1} - \eta^{n+1}_{m}
            + \eta^n_{m+1} - \eta^n_{m}
        \right] = 0,
    \end{array}
\end{equation}
where $a = \tau/(2h)$ and $H_k = H(x_k)$.

On the uniform orthogonal mesh
\[
\def\arraystretch{1.25}
\begin{array}{c}
    (x_0, \dots, x_M)\times(t_0, \dots, t_N),
    \\
    x_{k+1} - x_k = h = \textrm{const}, \qquad k = 0,\dots,M-1,
    \\
    t_{l+1} - t_l = \tau = \textrm{const}, \qquad l = 0,\dots,N-1,
\end{array}
\]
according to~(\ref{scm_sym_index}), we choose the following iterative process
\begin{equation} \label{scm_arb_iter}
    \def\arraystretch{2}
    \begin{array}{c}
            \overset{(j+1)}{\eta_m} =
            \eta^n_m - a \big[
                 \eta^{n}_{m} u^{n}_{m} + \overset{(j)}{\eta}_{m} \overset{(j)}{u}_{m}
                - \eta^{n}_{m-1} u^{n}_{m-1} - \eta^{n+1}_{m-1} \overset{(j)}{u}_{m-1}
                \\
                + (\overset{(j)}{u}_{m} + u^n_{m}) H_{m}
                - (\overset{(j)}{u}_{m-1} + u^n_{m-1}) H_{m-1}
            \big],
            \\
            \overset{(j+1)}{u_{m-1}} =
                u^n_{m-1}
                - a (
                    u_{m}^n \overset{(j)}{u}_{m}
                    - u_{m-1}^n \overset{(j)}{u}_{m-1}
                    + \eta^{n+1}_{m}
                    + \eta^{n}_{m}
                    - \eta^{n+1}_{m-1}
                    - \eta^{n}_{m-1}
                ),
            \\
            m = 1,\dots,M,
            \qquad
            n=0,\dots,N-1,
            \qquad
            j = 0, 1, 2, \dots,
    \end{array}
\end{equation}
where $H_k = H(x_0 + k h)$.
\noindent
Here for the $p^{\text{th}}$ time layer ($p=1,2,\dots$) we put
\[
    \overset{(0)}{u}_k = u^{p-1}_k,
    \qquad
    \overset{(0)}{\eta}_k = \eta^{p-1}_k,
    \qquad
    k = 0, \dots, M,
\]
and the values
\[
    \eta^p_0, u^p_M, \qquad p=0,\dots,N
\]
are determined by the initial and boundary conditions.

The iterative process is continued until
\[
\max\left\{
        \underset{k}{\max} \Big|\overset{(j+1)}{u}_{\!\!\!k}-\overset{(j)}{u}_k \Big|,
        \underset{k}{\max} \Big|\overset{(j+1)}{\eta}_{\!\!\!k}-\overset{(j)}{\eta}_k \Big|
    \right\}
    < \epsilon,
\]
for some $j \geqslant 0$
and some fixed $\epsilon$ ($0 < \epsilon \ll 1$).
Further on we choose~$\epsilon=10^{-6}$ for all the problems.

We also consider the viscous versions of the schemes, where the values in the process~(\ref{scm_arb_iter})
are modified as follows
\begin{equation}
 \overset{(j+1)}{\eta_m} \mapsto \overset{(j+1)}{\eta_m} - \nu\tau\!\dhm{D}(\eta_m),
 \qquad
 \overset{(j+1)}{u_{m}}
    \mapsto
 \overset{(j+1)}{u_{m}} - \nu\tau\!\dhm{D}(\eta_m u_m),
 \qquad
 m > 2,
\end{equation}
and $\nu$ is a linear artificial viscosity coefficient.
Introducing artificial viscosity often allows one to decrease
oscillations of numerical solutions.
For brevity, schemes without artificial viscosity
will also be called inviscid schemes.

\medskip

Further we state the problems for the schemes
on the river segment $x \in [0,L]$, where $L=100$.
The interval $L$ is uniformly divided into subintervals by the space step~$h=0.1$.
We choose the time step $\tau=0.1h=0.01$.

The parabolic bottom is defined as
\begin{equation} \label{parab_bottom_H}
  H(x) = d_1\left(\frac{2}{L}\right)^2\left(x - \frac{L}{2}\right)^2,
\end{equation}
where~$d_1$ is the depth of the fluid between the center of the parabola and some chosen zero level. We choose~$d_1 = 10$~throughout.

\medskip

\begin{figure}[ht]
    \centering
    \includegraphics[width=0.45\linewidth]{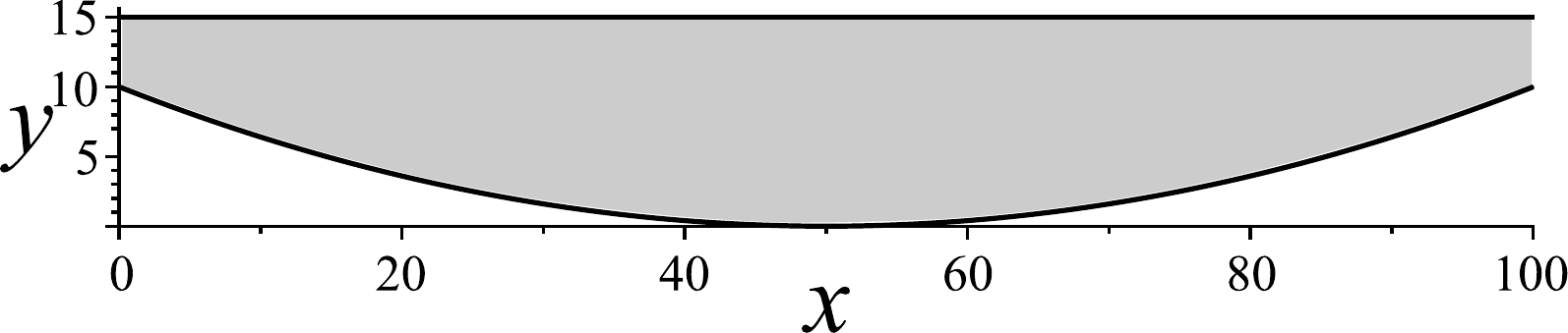}
    \caption{Parabolic bottom, a stationary solution~(no fluid flows).}
    \label{fig:parab-noflow}
\end{figure}

We consider, as the starting point, a stationary solution of the scheme~(\ref{scm_arb_iter})
for zero initial velocities of the fluid particles.
It is depicted in Figure~\ref{fig:parab-noflow}.
The height of the free surface above the chosen zero level is~$\eta=5$.
The conservation law of energy is held at any point,
so its plot is trivial and we do not provide it here.

\begin{figure}[ht]
    \centering
    \includegraphics[width=0.45\linewidth]{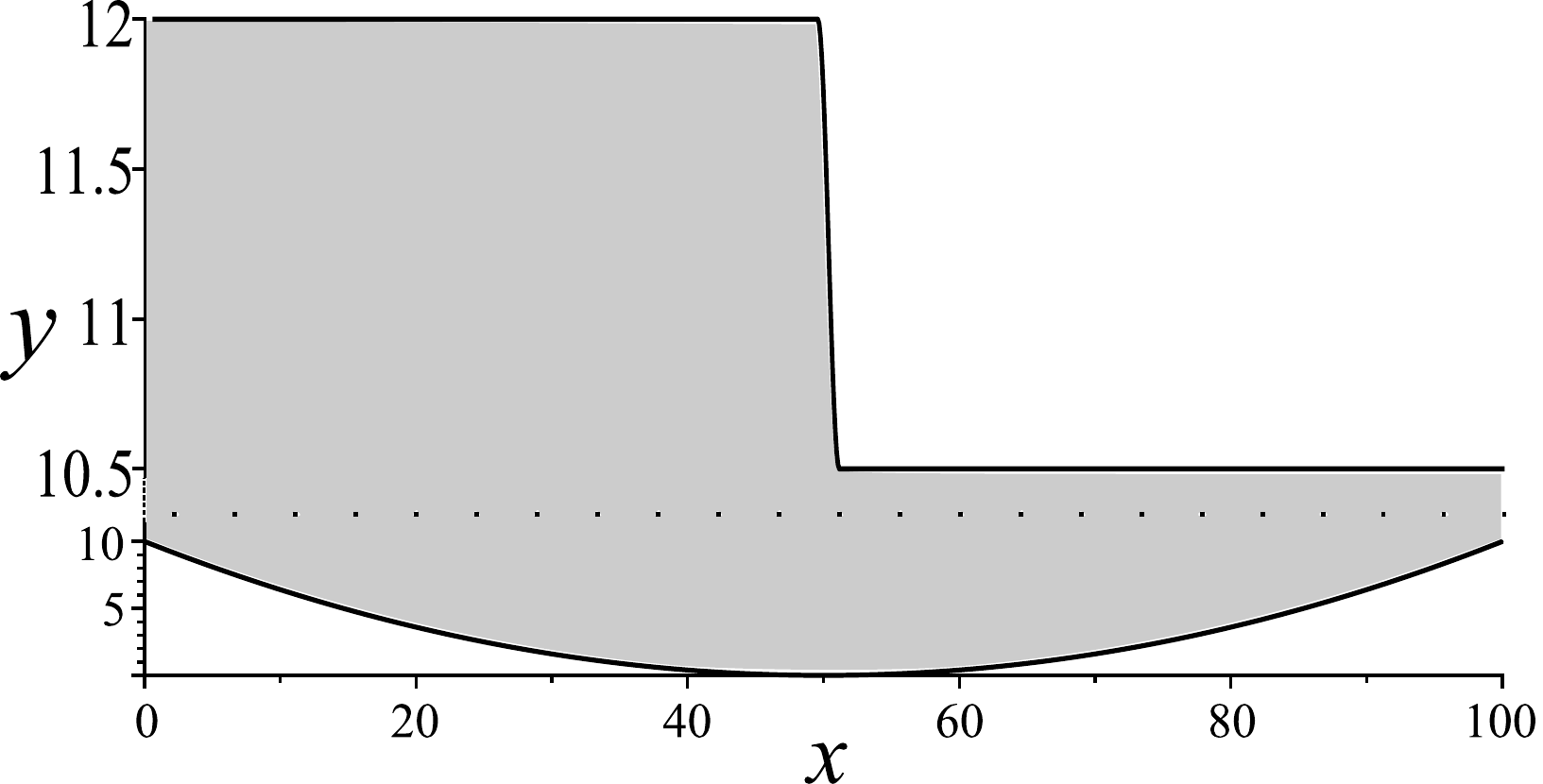}
    \caption{Parabolic bottom. The initial state~($t=0$) of the dam-break problem.}
    \label{fig:parab-dam-initial}
\end{figure}

\begin{figure}[ht]
    \centering
    \includegraphics[width=\linewidth]{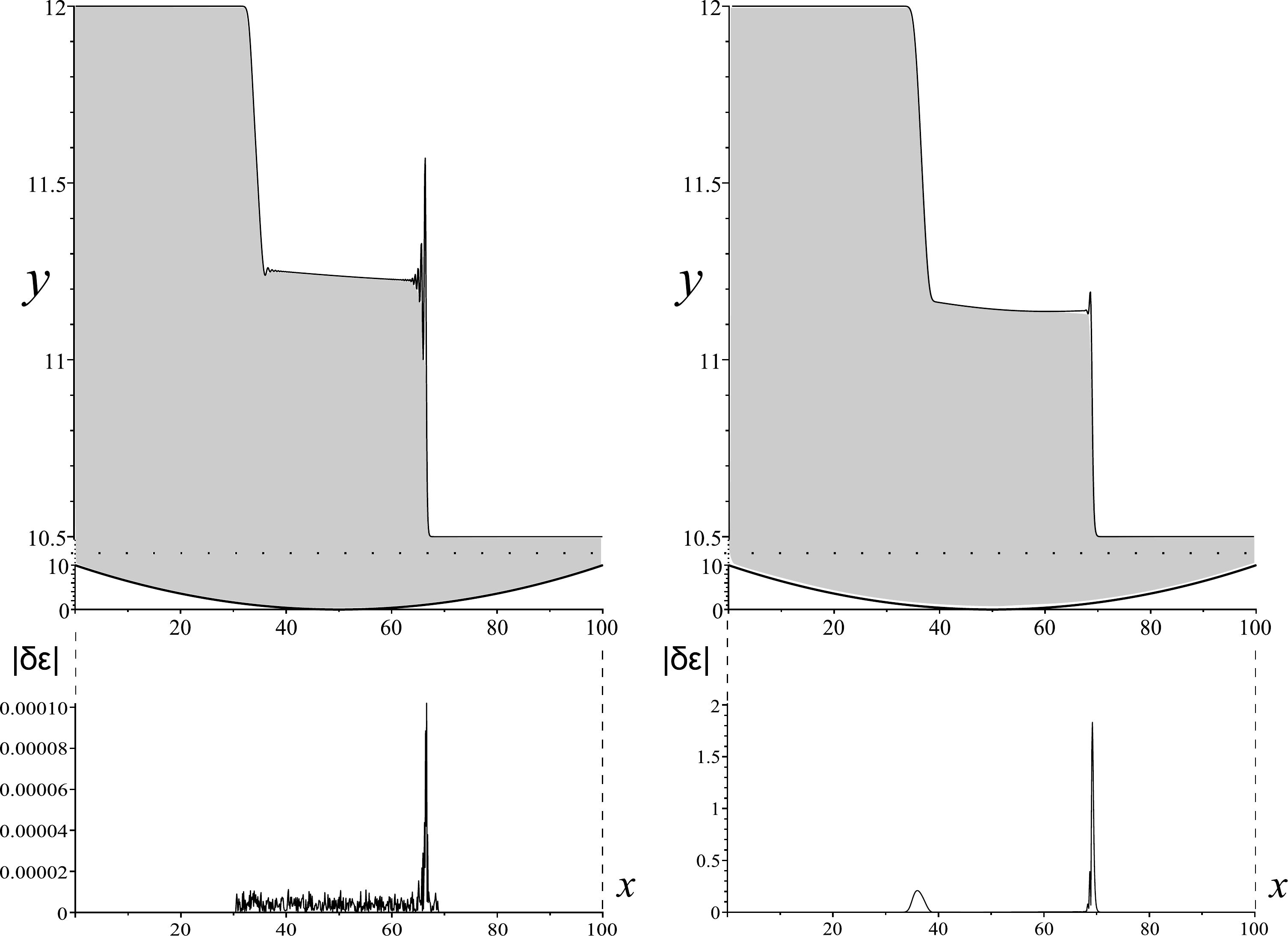}
    \caption{Parabolic bottom.
        Solution of the dam-break problem at $t=5.0$.
        \\
        Left: the inviscid version of the scheme.
        Right: the scheme with artificial viscosity~($\nu=0.08$).
        \\
        Free surface profiles are depicted on the top,
        and the conservation
        law of energy control values~$|\delta\varepsilon|$
        are depicted at the bottom of the figure.}
    \label{fig:parab-dam}
\end{figure}

\begin{figure}[ht]
\centering
\includegraphics[width=0.5\linewidth]{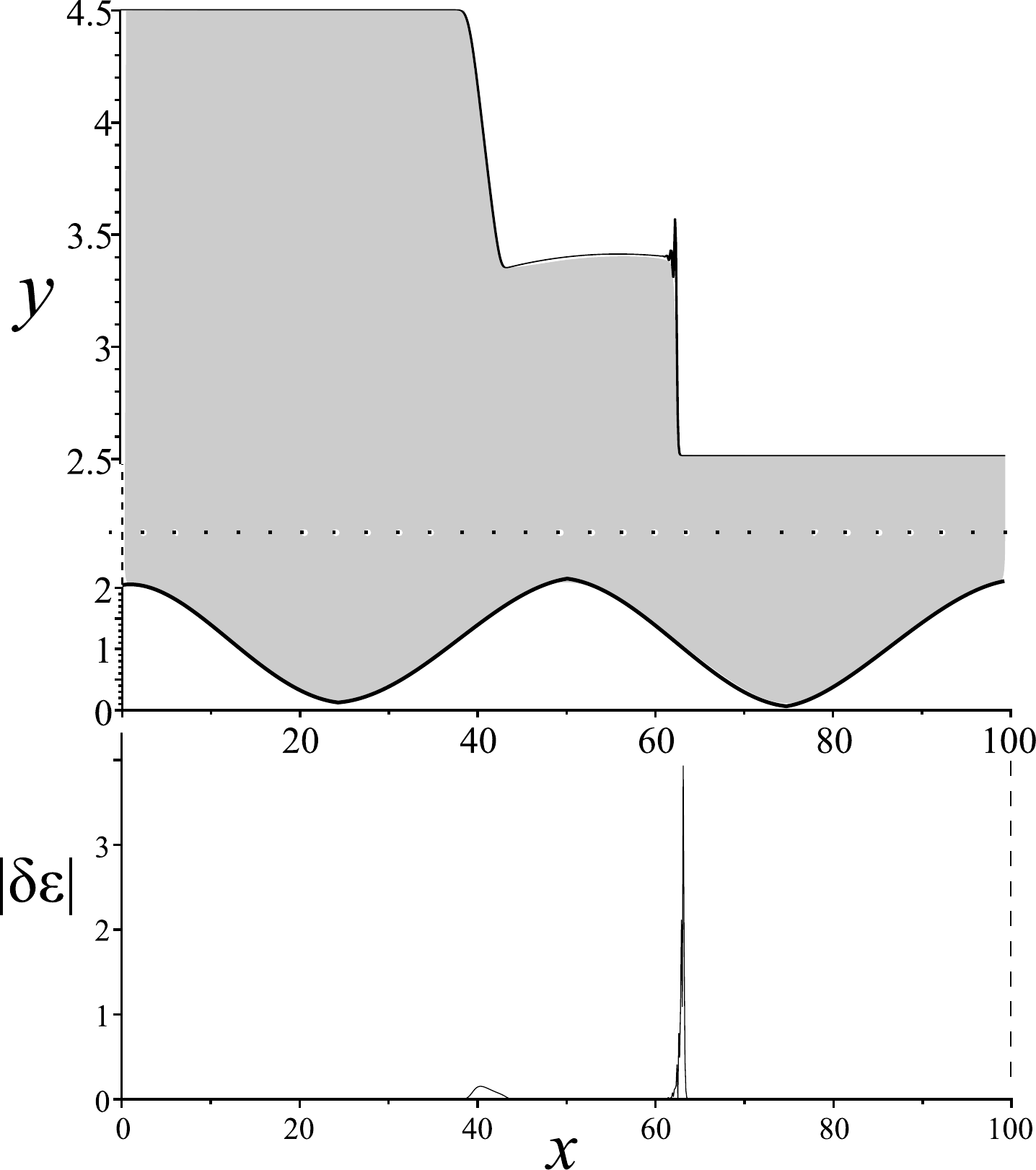}
\caption{Sinusoidal bottom.
                Solution of the dam-break problem at $t=5.0$ for
                a scheme with artificial viscosity~($\nu=0.15$).
                Free surface profile is depicted on the top,
                and the conservation law of energy control value~$|\delta\varepsilon|$
                is depicted at the bottom of the figure.}
\label{fig:sin-dam}
\end{figure}

Recall that~$\rho = \eta + H$ is the depth of the fluid above the bottom;
in all the pictures it is filled with light gray colour.
Notice also that the plots in Figures~\ref{fig:parab-dam-initial}--\ref{fig:sin-dam}
are scaled along the vertical axis. This allows one to see the free surface
profiles in more detail.

\medskip

Next we consider the dam-break problem for the parabolic bottom shape~(see Figure~\ref{fig:parab-dam-initial}).
The dam is located in the center~$x = L/2$ of the river segment.
At the initial moment of time the free surface heights to the left and to the right
of the dam are~$\eta_L = 2$ and~$\eta_R = 0.5$ appropriately.
Notice that the border between left and right segments
has been intentionally smoothed~(between 8 points).
We consider \emph{local} conservation laws, in order to avoid jump discontinuities
in the initial conditions.

In Figure~\ref{fig:parab-dam}, the solutions obtained
by the iterative processes~(\ref{scm_arb_iter})
are depicted for the schemes with and without artificial viscosity for~$t=5$.
The conservation of energy control values~$|\delta\varepsilon|$ are calculated as
absolute values of deviations of the finite-difference conservation law of energy 
from the zero value on actual solutions of the scheme.
We see from Figure~\ref{fig:parab-dam} that~$|\delta\varepsilon|$ takes on much greater
values for the inviscid scheme.

Finally, we consider the dam-break problem for a sinusoidal bottom shape, namely
\begin{equation}
H(x) = d_2 \cos^2\!\left(2\pi x/L\right),
\end{equation}
where we choose $d_2=2.0$, $\eta_L = 2.5$ and~$\eta_R = 0.5$.
The solution of the problem and the deviation~$|\delta\varepsilon|$
at $t=5.0$ for a viscous version of the scheme is given in Figure~\ref{fig:sin-dam}.

\subsection{Comparison with a modified non-conservative scheme}

In this section, we slightly modify scheme~(\ref{scm_sym})
so that it no longer possesses one of the conservation laws,
namely the conservation law of energy.
For this purpose, we change some constant
coefficients as follows~(the modified terms are underlined)
\begin{equation}\label{scm_naive}
    \def\arraystretch{1.75}
    \begin{array}{c}
        \dtp{D}(\eta_+)
        + \frac{1}{2}\dhp{D}\left(
            \eta u + \hat{\eta}\hat{u} + (\hat{u} + u) H
        \right) = 0,
        \\
        \dtp{D}(u) + \frac{1}{2} \dhp{D}\left(
            u\hat{u}
            + \underline{\tfrac{1}{2}\hat{\eta} + \tfrac{3}{2}\eta}
        \,\right) = 0.
    \end{array}
\end{equation}
Consequently, the conservation law of energy~(\ref{scm_sym_arb_cl}) transforms into
the following non-divergent equation
\begin{multline} \label{scm_naive_clE}
  \frac{1}{2}
    \dtp{D}\left(
        u^2 (\eta + H) + \eta_+^2
    \right)
    + \frac{1}{4}\dhp{D}\Big\{
    (u\hat{u} + \hat{\eta} + \eta)
    \left(
        \hat{u}\hat{\eta} + u\eta
        + (\hat{u} + u) H
    \right)
        + 2hu\hat{u}\!\dtp{D}(\eta)
    \Big\} =
    \\
   = \frac{\tau}{4} \left(\hat{u}\hat{\eta} + u\eta + (\hat{u} + u) H\right) \eta_{tx}.
\end{multline}
According to the results obtained in Section~\ref{sec:scheme_Euler_arb},
scheme~(\ref{scm_naive}) does not possess a local polynomial conservation law of energy.
It seems natural to consider~(\ref{scm_naive_clE}) as a
conservation law of energy approximation for scheme~(\ref{scm_naive}).
The right-hand side of~(\ref{scm_naive_clE}) essentially contributes to the energy dissipation of the scheme,
especially near high gradients of~$\eta$.
Notice that its terms cannot be represented as divergent expressions.
One can check this buy applying the Euler variational operator to the right-hand side of~(\ref{scm_naive_clE})
which should be zero in case of a divergent term.

\begin{figure}[ht]
\centering
\begin{minipage}[c]{0.45\textwidth}
\includegraphics[width=\linewidth]{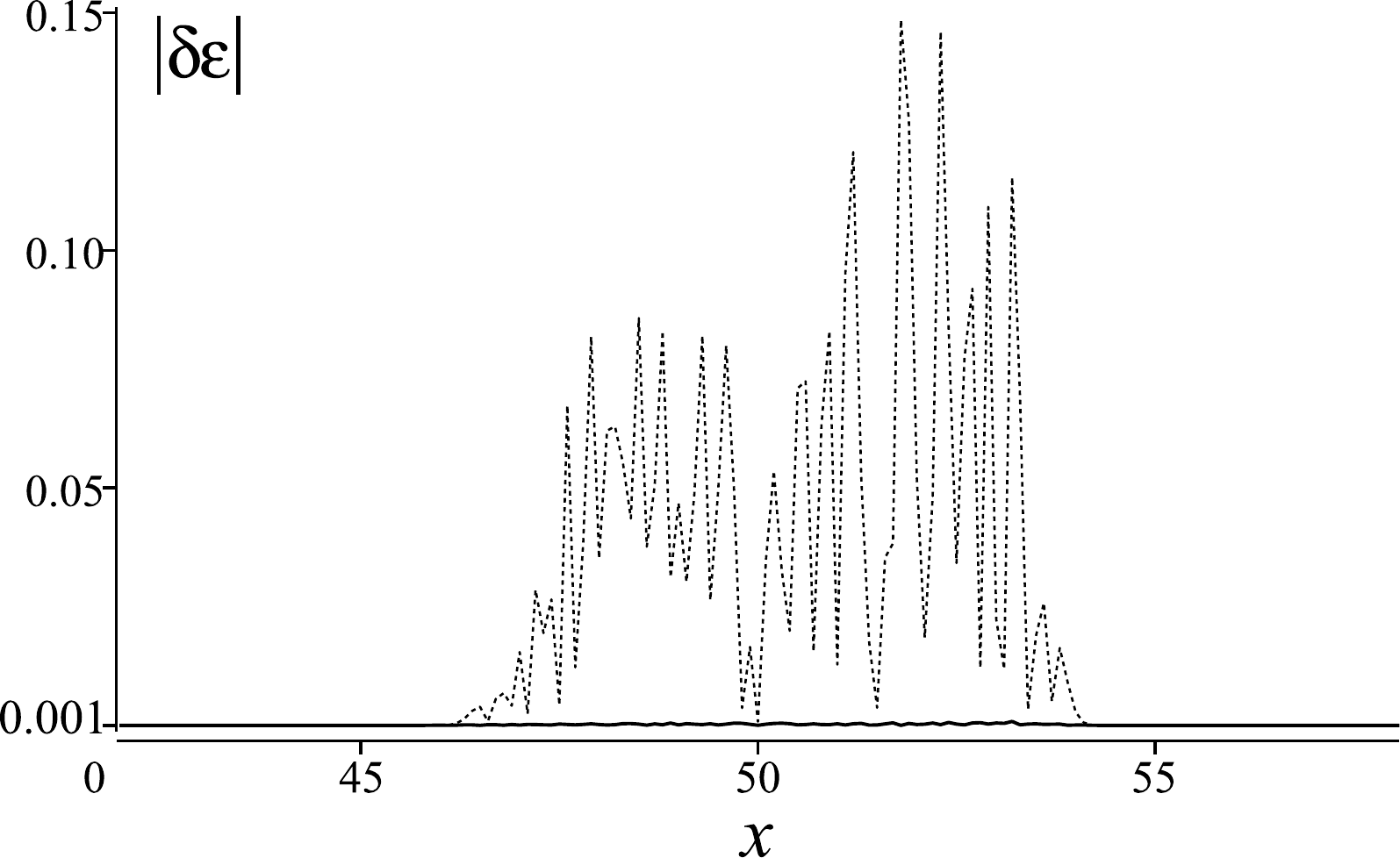}
\end{minipage}
\hfill
\begin{minipage}[c]{0.47\textwidth}
\includegraphics[width=\linewidth]{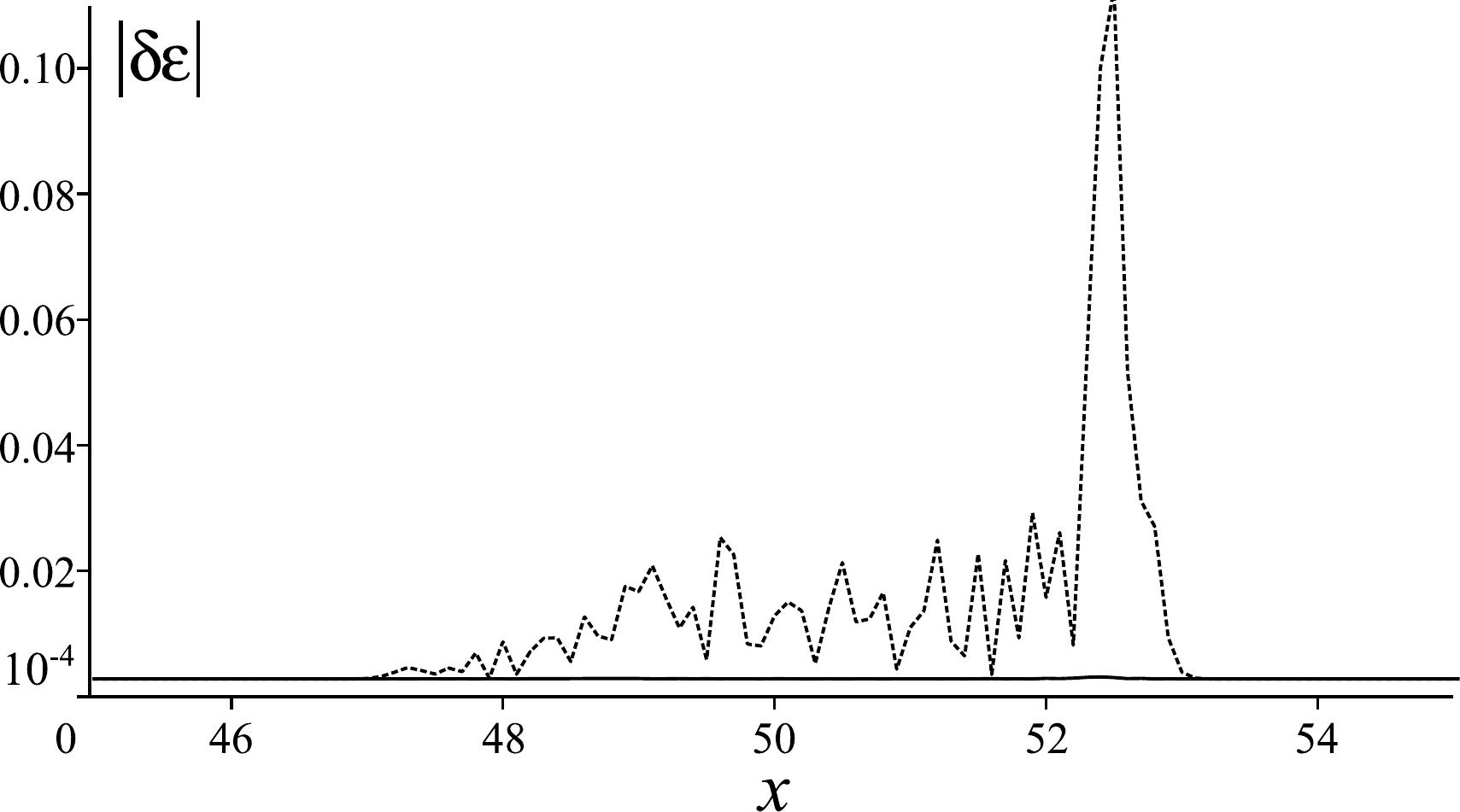}
\end{minipage}
\caption{Two inviscid schemes comparison (dam-break problem at $t=1.0$)
for parabolic (left) and sinusoidal (right) bottom shapes.
Solid lines~(---) correspond to the values of~$|\delta\varepsilon|$ for scheme~(\ref{scm_sym});
dot lines~($\cdots$) correspond to scheme~(\ref{scm_naive}).}
\label{fig:scm-cmd}
\end{figure}

In Figure~\ref{fig:scm-cmd} the values of~$|\delta\varepsilon|$ are given
for inviscid versions of schemes~(\ref{scm_sym}) and~(\ref{scm_naive})
for the example of the dam-break problem at~$t=1.0$.
Evidently, the deviation~$|\delta\varepsilon|$ for scheme~(\ref{scm_naive}) significantly
exceeds the deviation value for scheme~(\ref{scm_sym}).

\medskip

In addition we measure the change
of the total energy with time without regard to a particular scheme.
For this purpose, in accordance with~(\ref{Eul1}), we consider the sum
\begin{equation} \label{totalEapprox}
    \displaystyle
    \mathcal{H}(n) = \frac{h}{2}\sum_{(i)}\left[\rho^n_i(u^n_i)^2 + (\eta_i^n)^2\right]
\end{equation}
whose value on solutions of the schemes
should tend to constant in the continuous limit~\cite{bk:Dorodnitsyn[2011]}.
Recall that the value of~(\ref{totalEapprox}) corresponds to
the total energy for the shallow water equations that
is given\footnote{Actually, on a fixed segment $a \leqslant x\leqslant b$ the value may differ up to a constant which depends on $H$.} by the Hamiltonian~\cite{dorodnitsyn2019shallow, bk:Bihlo_numeric[2012]}
\begin{equation}
  \widetilde{\mathcal{H}} = \frac{1}{2}\int (\rho u^2 + \eta^2)dx.
\end{equation}

\begin{figure}[ht]
\centering
\begin{minipage}[b]{0.46\linewidth}
\includegraphics[width=\linewidth]{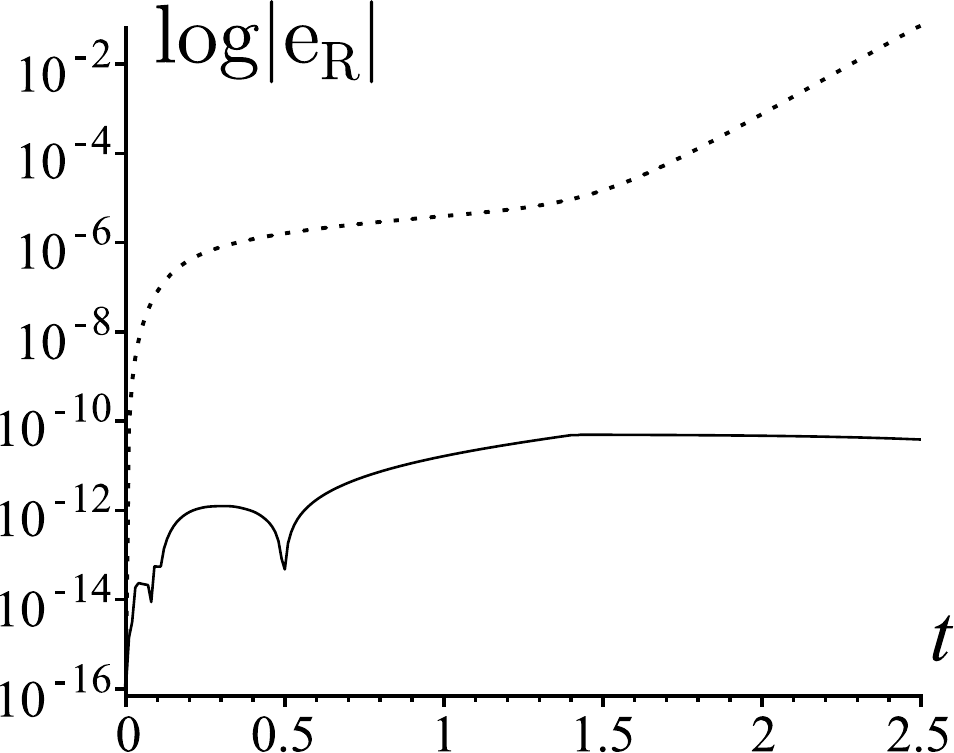}
\end{minipage}
\hfill
\begin{minipage}[b]{0.46\linewidth}
\includegraphics[width=\linewidth]{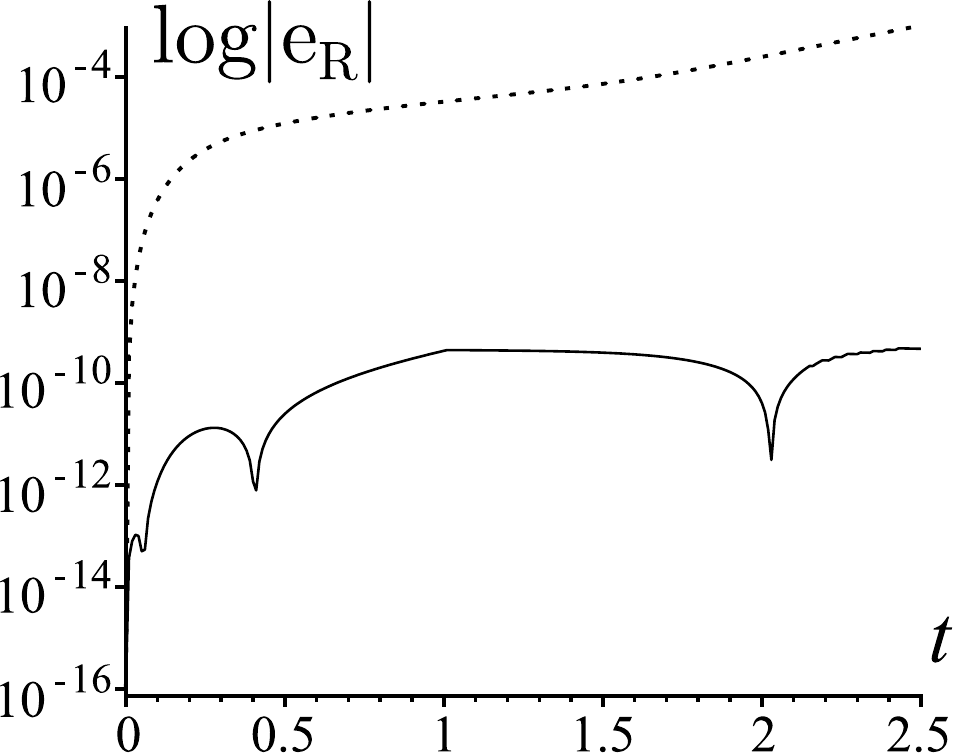}
\end{minipage}
\caption{Two inviscid schemes comparison (dam-break problem at $t=2.5$)
for parabolic (left) and sinusoidal (right) bottom shapes.
Solid lines~(---) correspond to the values of~$\log|e_R|$ for scheme~(\ref{scm_sym});
dot lines~($\cdots$) correspond to scheme~(\ref{scm_naive}).
The machine precision is~$10^{-16}$.}
\label{fig:scm-cmd2}
\end{figure}

The \emph{relative} change in energy is defined as follows
\begin{equation} \label{eRformula}
    e_R(n) = \frac{|\mathcal{H}(n) - \mathcal{H}(0)|}{|\mathcal{H}(0)|}.
\end{equation}
Its values~($0 \leqslant  t \leqslant 2.5$)
for the dam-break problem are given in Figure~\ref{fig:scm-cmd2}.
We see that scheme~(\ref{scm_sym}) conserves energy up to ten orders
of magnitude better than scheme~(\ref{scm_naive}).
It is remarkable that by adjusting the coefficients of the underlined
terms in~(\ref{scm_naive}) one can significantly improve the profile for~$\eta$
even without energy conservation, but the value of~$e_R$
still remains relatively high.

\section{Numerical implementation of invariant
conservative schemes for parabolic bottom topography
in Lagrangian coordinates}

The following scheme for a parabolic bottom in Lagrangian
coordinates is considered
\begin{equation} \label{scheme_sq_lagr_mod}
  \displaystyle
  \def\arraystretch{1.75}
  \begin{array}{c}
        \displaystyle
        x_{t\check{t}}
      + \frac{1}{2}\dsm{D}\left(
        \frac{1}{\hat{x}_s \check{x}_s}
      \right) + a_1 H^\prime(x)
      = 0,
      \\
      \displaystyle
      \qquad
      \tau_+ = \tau_-, \qquad
      h^s_+ = h^s_-,
  \end{array}
\end{equation}
where 
\[
H(x) = d_1\left[ \left(\frac{2}{L}\right)^2 \left(x - \frac{L}{2}\right)^2 - 1\right],
\]
and, according to Remark~\ref{rem:parabHequiv},
\[
a_1 = \frac{2 (\cosh(\sqrt{\beta_1}\tau) - 1)}{\beta_1\tau^2},
\qquad
\beta_1 = 8 d_1/L^2.
\]
Here $L$ is the length of the river segment and $d_1$ is the height of the parabolic bottom at the point~$x = L/2$.

Scheme~(\ref{scheme_sq_lagr_mod}) can be obtained from scheme~(\ref{scheme_sq_lagr})
by means of the following equivalence transformations~(see~\cite{bk:KaptsovMeleshko_1D_classf[2018]})
\begin{equation}
  x \mapsto \varepsilon_1 \varepsilon_2 x - \varepsilon_4,
  \qquad
  t \mapsto \varepsilon_3 t,
  \qquad
  s \mapsto \varepsilon_3 s,
  \qquad
  H^\prime \mapsto -\varepsilon_2 H^\prime / (\varepsilon_3)^{2},
\end{equation}
where
\[
    \varepsilon_1 = \frac{2^{1/3} L^2}{8 d_1},
    \qquad
    \varepsilon_2 = \frac{16\sqrt{2} d_1^{3/2}}{L^3},
    \qquad
    \varepsilon_3 = \frac{\sqrt{8 d_1}}{L},
    \qquad
    \varepsilon_4 = -2^{5/6} \sqrt{d_1}.
\]
Notice that for scheme~(\ref{scheme_sq_lagr_mod}) the conservation low of energy~(\ref{CLparabEnergy})
has the following form
\begin{equation}
    \dtm{D}\left[
        x_t^2 + \frac{1}{2}(x_s^{-1} + \hat{x}_s^{-1})
        + a_1\beta_1\left(x - \frac{L}{2}\right)\left(\hat{x} - \frac{L}{2}\right)
    \right]
    + \dsm{D}\left(
        \frac{x_t^+ + \check{x}_t^+}{2\hat{x}_s \check{x}_s}
    \right) = 0.
\end{equation}

\medskip

The statement of initial and boundary conditions is related to the solution
of the following Cauchy problem~\cite{bk:YanenkRojd[1968]}
\begin{equation}
x_t(t,\xi) = u(t, x(t, \xi)),
\qquad
x(t_0, \xi) = \xi,
\end{equation}
where $\xi = x|_{t = t_0}$ is Lagrangian label of a particle.
The conservation law of mass in Lagrangian coordinates is
\begin{equation}
\rho(t, x(t, \xi)) = \frac{\rho_0(\xi)}{x_\xi(t, \xi)},
\end{equation}
where $\rho_0(\xi)$ is an arbitatry function of integration.
The mass Lagrangian coordinate~$s$ is introduced by the relation~$\xi = \alpha(s)$
in such a way that
\begin{equation}
\tilde{\rho}(t, s) = \rho(t, x(t, \alpha(s))) = \frac{1}{\tilde{x}_s(t,s)}.
\end{equation}
The sign $\tilde{\,}$ is further omitted.

Thus, given an initial height~$\rho(\xi)$ of the fluid over the bottom~$H$,
one obtains the function~$\alpha(s)$ by solving the Cauchy problem
\begin{equation} \label{CauchyProblem}
\rho(\alpha(\xi)) = \frac{1}{\alpha^\prime(s)},
\qquad
\alpha(0) = 0.
\end{equation} \label{SofXintegral}
The initial distribution~$s(x)$ is obtained either
as inverse~$s = \alpha^{-1}(x)$
or directly by integrating the equation
\begin{equation}
s(x)|_{t=t_0} = \int_0^x \rho(\xi) d\xi.
\end{equation}
Numerical computations in mass Lagrangian coordinates are
considered in the paper~\cite{bk:DorKapMelGN2020}
for the example of the Green-Naghdi equations in more detail.

\medskip

In order to implement scheme~(\ref{scheme_sq_lagr_mod}) numerically,
we represent its first equation in the following form
\begin{equation} \label{scheme_parab_nonlin_form}
\hat{x} - 2x + \check{x}
    + \frac{h\tau^2}{2} \frac{
        (\hat{x} - \hat{x}_-)(\check{x} - \check{x}_-)
        -
        (\hat{x}_+ - \hat{x})(\check{x}_+ - \check{x})
    }
    {
        (\hat{x} - \hat{x}_-)(\check{x} - \check{x}_-)
        (\hat{x}_+ - \hat{x})(\check{x}_+ - \check{x})
    }
    + \tau^2 a_1 H^\prime(x) = 0,
\end{equation}
Then, we linearize equation~(\ref{scheme_parab_nonlin_form}) by representing it as the following
iterative procedure
\begin{multline*}
    x^{(j+1)}_m - 2x^n_m + x^{n-1}_m
    + \tau^2 a_1 H^\prime(x^n_m)
    \\
    + \frac{h^s\tau^2}{\Delta_1} \left(
        (x^{(j+1)}_m - x^{(j+1)}_{m-1})(x^{n-1}_m - x^{n-1}_{m-1})
        -
        (x^{(j+1)}_{m+1} - \hat{x})(x^{n-1}_{m+1} - x^{n-1}_m)
    \right) = 0
\end{multline*}
or
\begin{multline} \label{tridiag}
    \frac{h^2\tau^2}{\Delta_1}(x^{n-1}_m - x^{n-1}_{m-1}) \, x^{(j+1)}_{m-1}
    - \left(1 + \frac{h^2\tau^2}{\Delta_1}(x^{n-1}_{m+1} - x^{n-1}_{m-1}) \right) x^{(j+1)}_{m}
    \\
    + \frac{h^2\tau^2}{\Delta_1}(x^{n-1}_{m+1} - x^{n-1}_m) \, x^{(j+1)}_{m+1}
    = 2x^n_m - x^{n-1}_m -\tau^2 a_1 H^\prime(x^n_m),
\end{multline}
where
\[
\Delta_1 = 2\,(x^{(j)}_m - x^{(j)}_{m-1})
        (x^{(j)}_{m+1} - x^{(j)}_m)
        (x^{n-1}_m - x^{n-1}_{m-1})
        (x^{n-1}_{m+1} - x^{n-1}_m),
\]
\[
    n = 2, 3, \dots,
    \qquad
    m = 2, 3, \dots, \lfloor L/h^s \rfloor - 1,
\]
and the indices ${}^{(j)}$ denote the number of iteration.
System~(\ref{tridiag}) on each iteration
can be solved with the help of tridiagonal matrix algorithm\footnote{See details on
this well-known method and its stability conditions, for example, in~\cite{bk:Samarskii2001theory}.
It is easy to verify that the stability conditions of the algorithm are satisfied for system~(\ref{tridiag}).}.
Notice that scheme~(\ref{scheme_sq_lagr_mod}) requires
two time layers
\begin{equation} \label{X0X1}
    x^0_m \quad \text{and} \quad x^1_m, \qquad m = 1, 2, \dots, \lfloor L/h^s \rfloor,
\end{equation}
on which the initial conditions must be specified.
In the following examples, the first and second layers
can be considered equal, which greatly simplifies
the setting of the initial conditions.

\subsection{A stationary flow over the parabolic bottom}

\begin{figure}[H]
\centering
\includegraphics[width=0.25\linewidth]{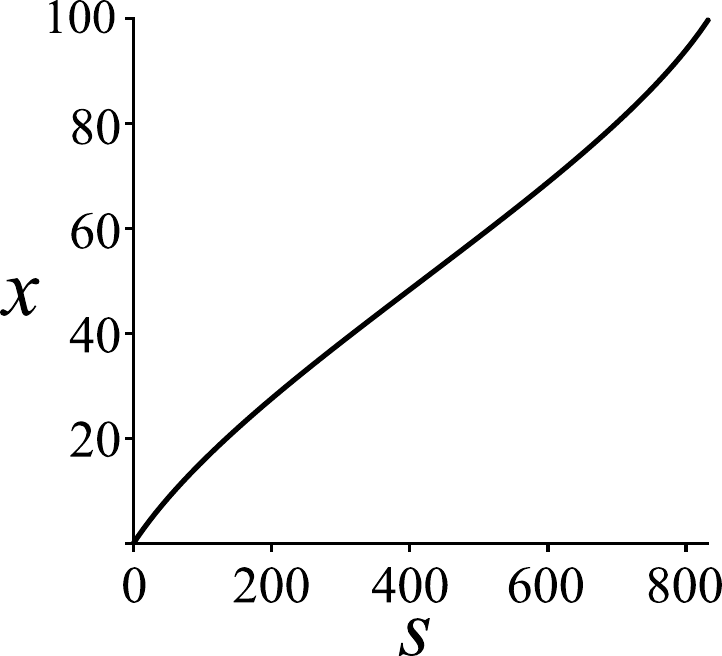}
\caption{The relation~$x(s)$ on the segment $0 \leqslant x \leqslant L$.}
\label{fig:s_of_x_flat}
\end{figure}

Consider a stationary solution of the scheme~(\ref{scheme_sq_lagr_mod})
for zero initial velocities of the fluid particles.
This problem was considered in Eulerian coordinates in the beginning of Section~\ref{sec:Euler_scm_impl}.
The initial height of the fluid over the bottom is
\begin{equation}
  \rho(\xi) = \eta - d_1\left(\frac{2}{L}\right)^2\left[\left(\xi - \frac{L}{2}\right)^2 - \left(\frac{L}{2}\right)^2 \right].
\end{equation}
The numerical solution of the Cauchy problem~(\ref{CauchyProblem}) is given in Figure~\ref{fig:s_of_x_flat}.
This allows one to specify the initial distribution of the fluid particles.

\begin{figure}[ht]
\centering
\includegraphics[width=0.75\linewidth]{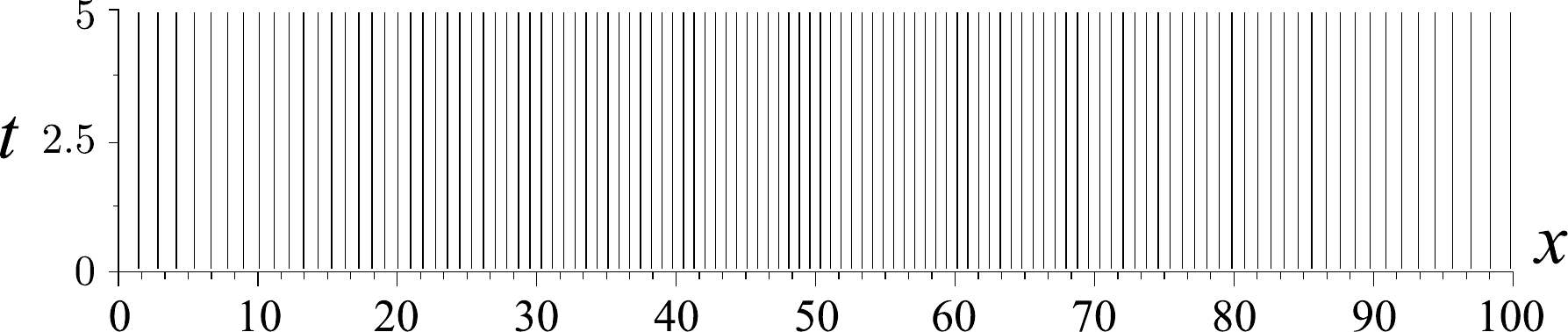}
\caption{The flow of the particles~$x(t)$ for~$0 \leqslant t \leqslant 5$.}
\label{fig:x_flow_flat}
\end{figure}

In case of zero initial velocities, the values of~$x$ on layers~(\ref{X0X1}) coincide,
and thus the initial conditions can easily be posed.
The stationary flow for $\eta=5$, $d_1=10$, $h^s=0.1$ and $L=100$ is depicted in Figure~\ref{fig:x_flow_flat}.
Notice that mass particles are located the denser to another the closer they are to the center of the parabolic bottom.

\subsection{The dam-break problem over the parabolic bottom}

\begin{figure}[ht]
\centering
\begin{minipage}[c]{0.30\linewidth}
    \includegraphics[width=\linewidth]{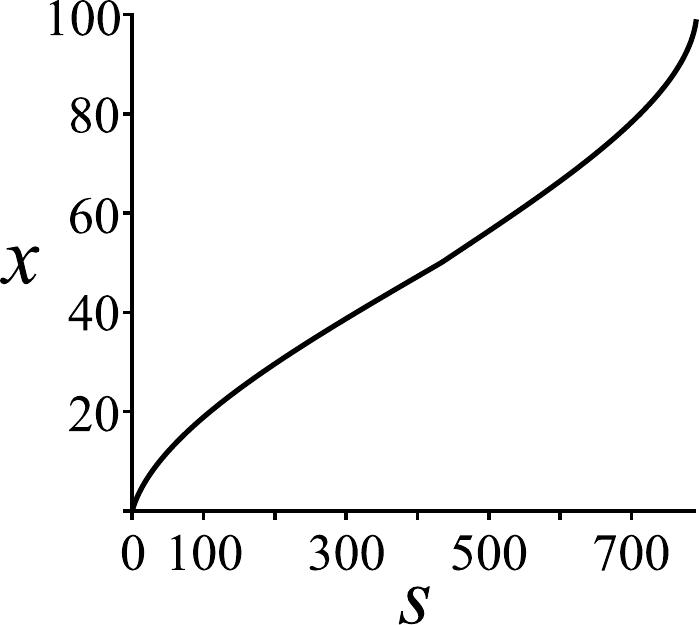}
  \end{minipage}
  \hfil
  \begin{minipage}[c]{0.58\linewidth}
  \includegraphics[width=\linewidth]{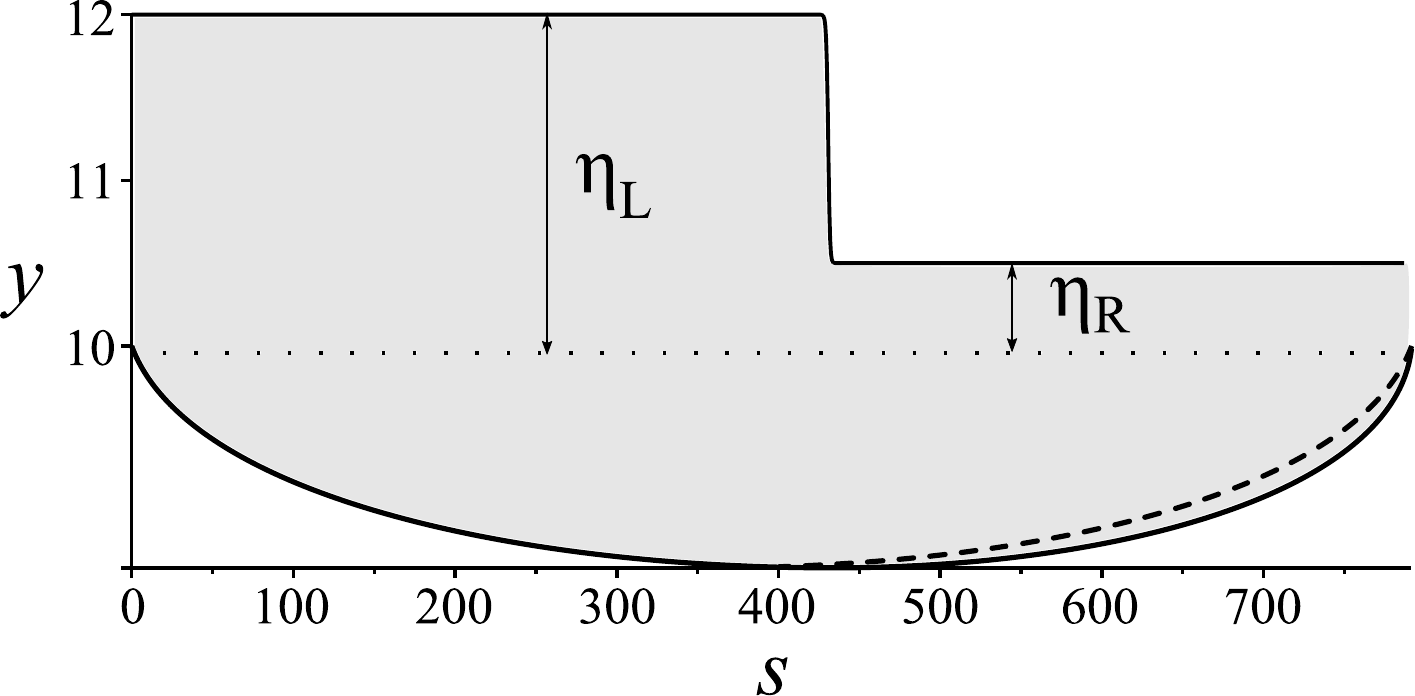}
  \end{minipage}
\caption{The relation~$x(s)$ on the segment $0 \leqslant x \leqslant L$ (left)
and the bottom profile~$\widetilde{H}(s) = H(x(s))$ in~Lagrangian coordinates~(right).
The dashed line demonstrates the asymmetry of the bottom profile in case of~Lagrangian coordinates.
This asymmetry occurs due to large discontinuity in the initial data.}
\label{fig:x_initial}
\end{figure}

Consider the dam-break problem which was stated in Section~\ref{sec:Euler_scm_impl}.
In order to integrate~(\ref{SofXintegral}) analytically and to
provide smoother initial data for the scheme, we approximate
the initial free surface profile by the function
\begin{equation}
\eta(\xi) = \eta_L - \frac{\eta_R - \eta_L}{1 + \exp\left(\sigma_1(L/2 - \xi)\right)},
\end{equation}
where $\sigma_1 = 20$ is the curve steepness coefficient, and the constants~$\eta_L=2$ and~$\eta_R=0.5$
are given in~Figure~\ref{fig:x_initial}~(right).
As in the pervious example, we put $d_1=10$ and $L=100$.

The solution of the corresponding Cauchy problem~(\ref{CauchyProblem}) is
given in Figure~\ref{fig:x_initial}~(left).
This solution corresponds to the layer~$x^0_m$ of~(\ref{X0X1}).
The second time layer has the same values as we
set the initial velocities to zero and considering
the behavior of the liquid under gravity.
Notice that in mass Lagrangian coordinates the bottom
function~$\widetilde{H}(s) = H(x(s))$
may differ from the function~$H(x)$ in Eulerian
coordinates~(see the right side of Figure~\ref{fig:x_initial}).

%

\begin{figure}[ht]
\centering
\includegraphics[width=0.75\linewidth]{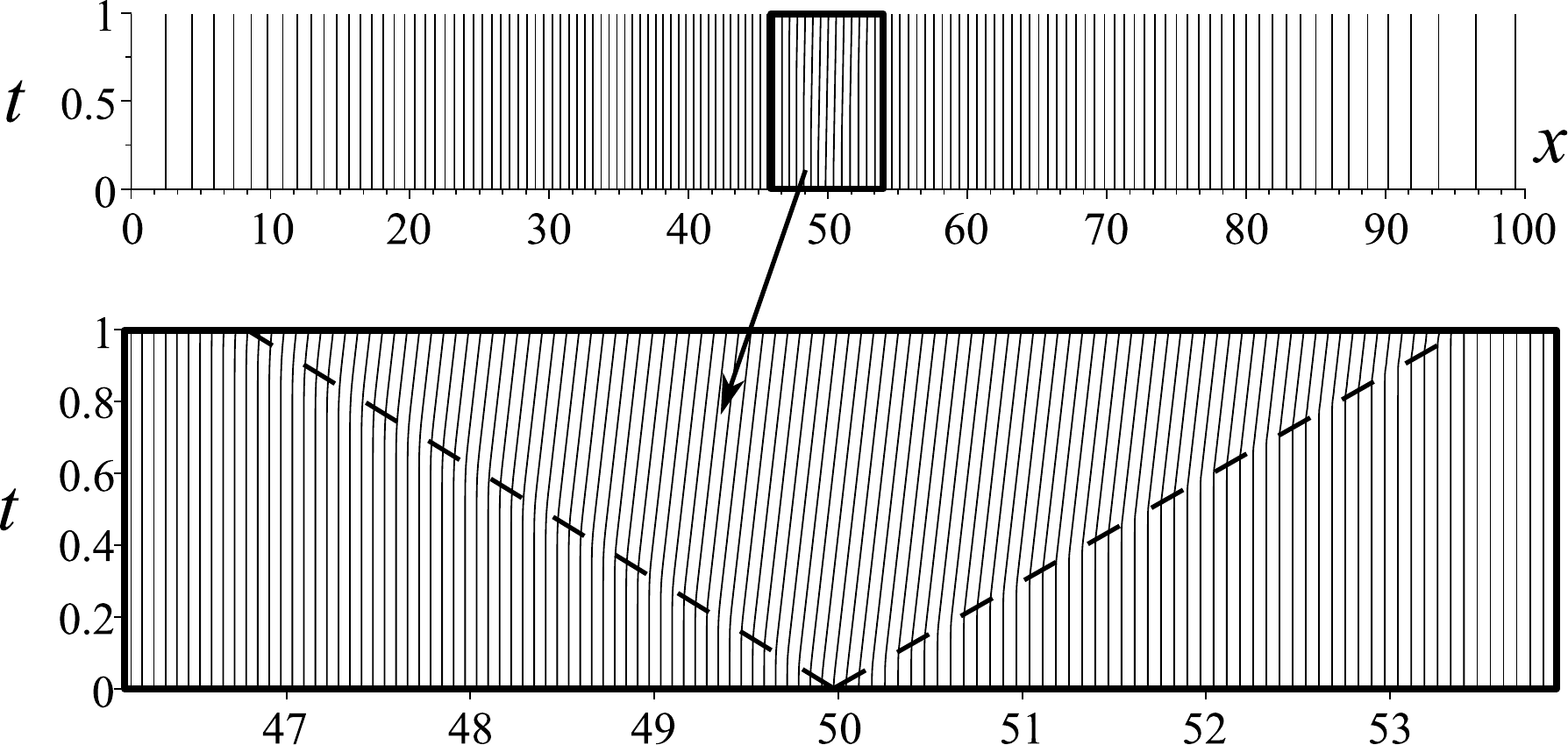}
\caption{The flow of the particles~$x(t)$ for the parabolic bottom~($0 \leqslant t \leqslant 1$).
The dashed lines outline the characteristics of the flow.}
\label{fig:x_flow_parab}
\end{figure}

The resulting flow of the fluid particles for $t=1$, $h^s = 0.25$ and~$\tau=0.00125$ is given in Figure~\ref{fig:x_flow_parab}.

Along with scheme~(\ref{scheme_sq_lagr_mod}) we consider its viscous version
which includes an artificial linear viscosity term~\cite{bk:SamarskyPopov_book[1992]}.
In hydrodynamic coordinates this term is proportional to~$\rho u_s$,
and in mass Lagrangian coordinates, by means of~(\ref{s}),
it should be proportional to~$x_{ts}/x_s$.
Various finite-difference representations for artificial viscosity term are possible
and we choose the following one
\begin{equation} \label{scheme_sq_lagr_mod_visc}
  \displaystyle
  \def\arraystretch{1.75}
  \begin{array}{c}
        \displaystyle
        x_{t\check{t}}
      + \frac{1}{2}\dsm{D}\left(
        \frac{1}{\hat{x}_s \check{x}_s}
      \right)
      + a_1 H^\prime(x)
      -\nu \frac{\check{x}_{ts}}{x_s}
      = 0,
      \\
      \displaystyle
      \qquad
      \tau_+ = \tau_-, \qquad
      h^s_+ = h^s_-,
  \end{array}
\end{equation}
where $\nu \sim h$ is some constant linear viscosity coefficient.

\begin{figure}[H]
\centering
\includegraphics[width=0.67\linewidth]{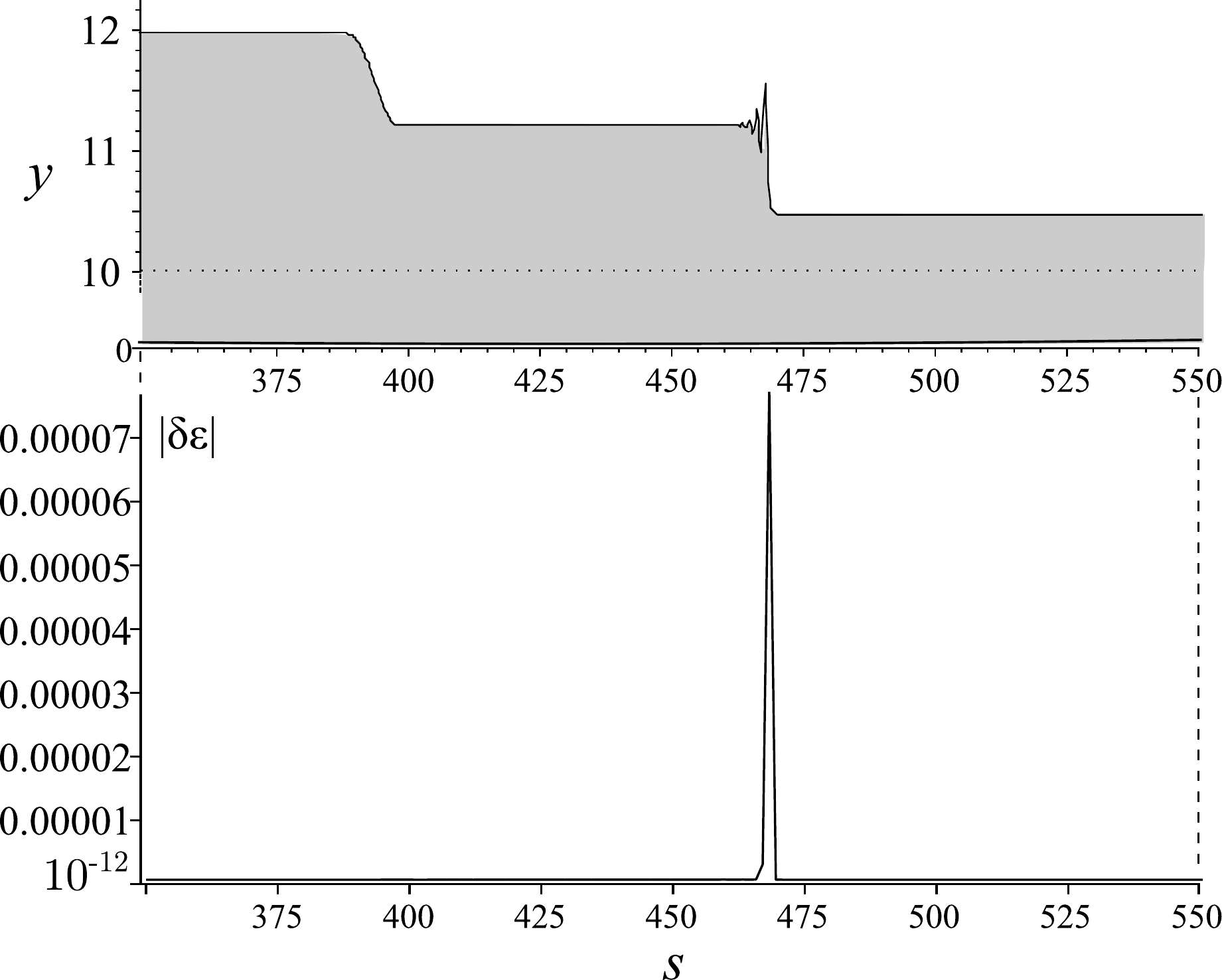}
\caption{Parabolic bottom.
        Solution $\rho=1/x_s$ of the dam-break problem at $t=1$ for the inviscid scheme
        on the segment~$350 \leqslant s \leqslant 550$.
        Free surface profile is depicted on the top,
        and the conservation
        law of energy control value~$|\delta\varepsilon|$~(in logarithmic scale)
        is depicted at the bottom of the figure.}
\label{fig:lagr-parab-rho}

\vspace*{\floatsep}

\includegraphics[width=0.66\linewidth]{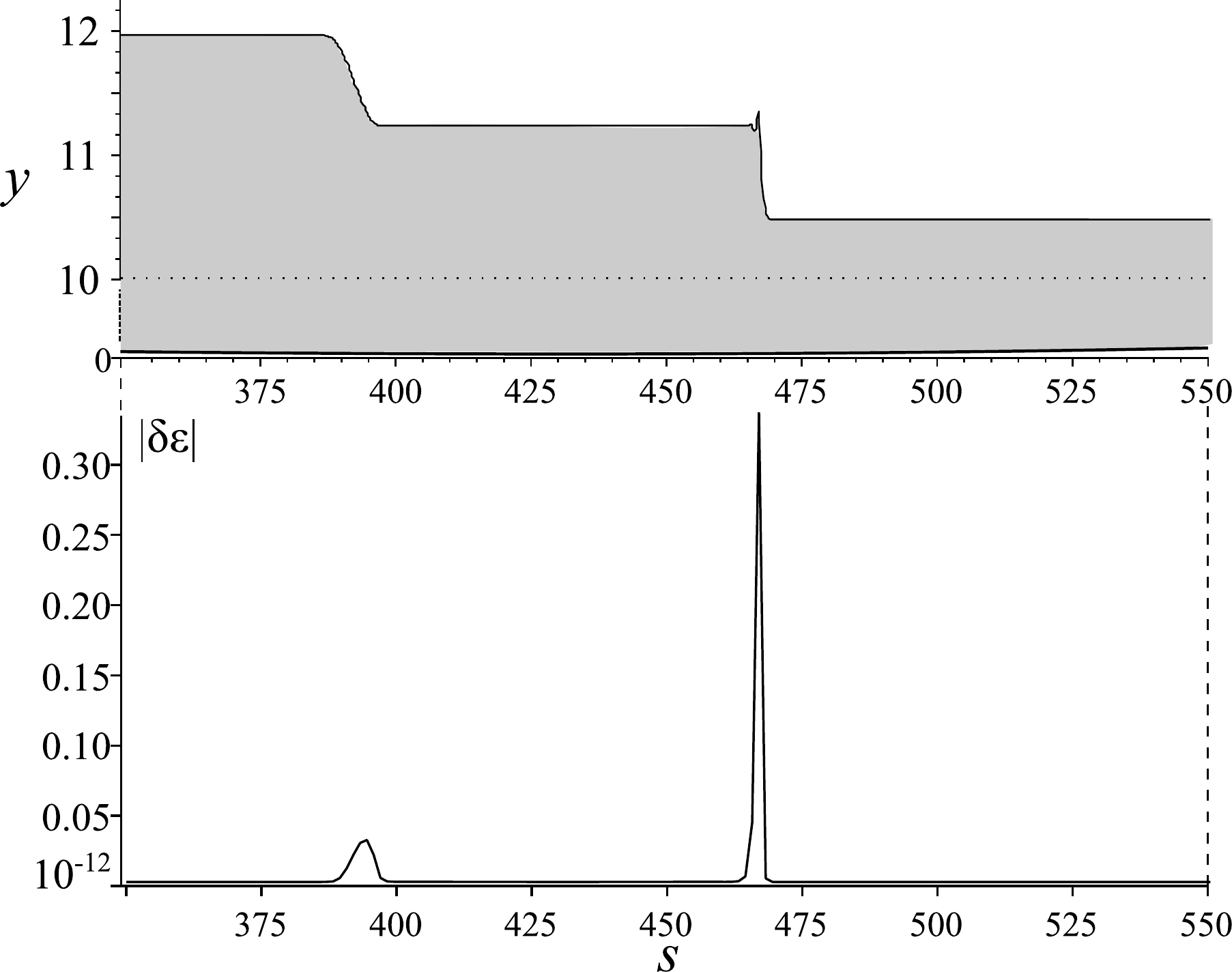}
\caption{In contrast with the results presented in Figure~\ref{fig:lagr-parab-rho},
    a scheme with artificial viscosity~($\nu = h$) allows one to get a better free surface profile.
    At the same time, the energy is not so well conserved.}
\label{fig:lagr-parab-rho-visc}
\end{figure}

The solutions of the dam-break problem at $t=1$ for schemes~(\ref{scheme_sq_lagr_mod})
and~(\ref{scheme_sq_lagr_mod_visc}) are presented
in~Figures~\ref{fig:lagr-parab-rho} and~\ref{fig:lagr-parab-rho-visc}.
The results of the numerical calculations are similar to the results obtained
in~Eulerian coordinates for the case of parabolic bottom topography.

\section{Conclusion}

Invariant finite-difference schemes for the shallow water equations in Lagrangian coordinates for the cases of flat bottom, inclined and arbitrary bottom topographies were considered in~\cite{dorodnitsyn2019shallow}.
In the present paper this research is extended to the schemes for parabolic bottom topography.
Invariant schemes which possess the conservation laws of mass, energy, momentum and center-of-mass law, as well as two additional conservation laws, are constructed on three time layers. It is shown that the schemes can be written on two time layers in terms of  hydrodynamic variables. A numerical implementation of one of the schemes in Lagrangian coordinates, as well as its version with pseudo-viscosity, is performed. The dam-break problem is considered as the main example. Calculations show that energy is well conserved on the obtained solutions, while the conservation of energy depends significantly on the pseudo-viscosity value.

The main part of the publication is devoted to the construction of difference schemes in~Eulerian
coordinates for an arbitrary bottom topography.
By means of a finite-difference analogue of the direct method,
an invariant conservative scheme is constructed on a uniform orthogonal mesh.
This scheme possesses difference analogues of all the differential conservation laws: mass, energy, momentum conservation and center-of-mass law. Using the example of the dam-break problem, a numerical analysis of the constructed scheme is carried out. It is shown that the constructed scheme has obvious advantages in conservativity:~even a small modification of the scheme can significantly worsen energy preservation.

On the basis of the classifications~\cite{bk:KaptsovMeleshko_1D_classf[2018],bk:AksenovDruzkov_classif[2019]}
symmetries and conservation laws in Eulerian and Lagrangian coordinates for the finite-difference
cases are presented in~Tables~\ref{tab:Htable} and~\ref{tab:HtableDelta}.
These tables cover both the schemes considered in the present paper and in the paper~\cite{dorodnitsyn2019shallow}.


\section*{Acknowledgements}

The research was supported by Russian Science Foundation Grant No 18-11-00238
``Hydro\-dyna\-mics-type equations: symmetries, conservation
laws, invariant difference schemes''.
E.K. acknowledges Suranaree University of Technology for
Full-time Master Researcher Fellowship.
The authors are grateful to S.~V.~Meleshko and E.~Schulz for valuable discussions.

\bibliographystyle{unsrt}

%

\begin{landscape}
\bgroup
\def\arraystretch{1.5}
\begin{longtable}[ht]{|c|c|c|c|c|}
\caption{Symmetries and conservation laws for differential systems in Eulerian coordinates}
\label{tab:Htable}\\
\hline
\# & Bottom $H$ &
    $\begin{array}{c}\textrm{Symmetries of differential}\\
    \textrm{system~(\ref{Euler1}),(\ref{Euler2})}\end{array}$
&
    $\begin{array}{c}\textrm{Related differential}\\ \textrm{conservation laws}\end{array}$
\\ \hline\hline
\endfirsthead
\multicolumn{3}{c}%
{{\it Table \thetable\ Continued from previous page}} \\
\hline
\# & Bottom $H$ &
    $\begin{array}{c}\textrm{Symmetries of differential}\\
    \textrm{system~(\ref{Euler1}),(\ref{Euler2})}\end{array}$
&
    $\begin{array}{c}\textrm{Related differential}\\ \textrm{conservation laws}\end{array}$
\\ \hline\hline
\endhead
-- &

    $\begin{array}{c}
        \textrm{arbitrary}\\
       H(x)
     \end{array}$
 &
 $X_1 = {\frac{\partial}{\partial{t}}}$
 &
 \begingroup
 \renewcommand*{\arraystretch}{1.5}
 \begin{tabular*}{0.52\textwidth}{l|l}
    & \mmomconservtxt\\ \hline
    $X_1$
    &
    $\begin{array}{l}
    \text{\energconservtxt:}\\
    D_t(u^2 (\eta + H) + \eta^2)
    \\
    \quad +D_x((\eta + H)(u^3 + 2 u \eta)) = 0
 \end{array}$
 \end{tabular*}
 \endgroup
 \\
 \hline
1 &
 $\begin{array}{c}
    0
  \end{array}$
 &
 $\begin{array}{c}
    X^1_2 = t {\frac{\partial}{\partial{t}}} + x {\frac{\partial}{\partial{x}}},
    \\
    X^1_3 = x {\frac{\partial}{\partial{x}}} + u {\frac{\partial}{\partial{u}}} + 2\eta {\frac{\partial}{\partial\eta}},
    \\
    X^1_4 = t {\frac{\partial}{\partial{x}}} + {\frac{\partial}{\partial{u}}},
    \\
    X^1_5 = \left(\frac{1}{2}x - \frac{3}{2}t u \right) {\frac{\partial}{\partial{t}}}
    + \left(\frac{3}{2} t\eta - \frac{3}{4}t u^2 \right) {\frac{\partial}{\partial{x}}}
    \\
    \qquad
    + \left(\frac{1}{4} u^2 + \eta \right) {\frac{\partial}{\partial{u}}}
    + u \eta {\frac{\partial}{\partial\eta}}
    \\
    X_\infty^1 = w_1 {\frac{\partial}{\partial{t}}} + (u w_1(u, \eta) + w_2(u, \eta)){\frac{\partial}{\partial{x}}}
    \\
    (X_\infty^1 \supset {\frac{\partial}{\partial{x}}})
 \end{array}$
 &
 \begingroup
 \renewcommand*{\arraystretch}{5}
 \begin{tabular*}{0.52\textwidth}{l|l}
    & \mmomconservtxt\\ \hline
   $X_1$ & \energconservtxt
    \\
    \hline
   ${\frac{\partial}{\partial{x}}}$
   &
   $D_t(u \eta) + D_x(u^2 \eta + \frac{1}{2} \eta^2) = 0$
   \\
 \end{tabular*}
 \endgroup
 \\
 \hline
 2 &
 $\begin{array}{c}
    x
 \end{array}$
 &
 \multicolumn{2}{c|}{Can be transformed to the flat bottom by~(\ref{bottom_tr_Chirkunov})}
 \\
 \hline
3 &
 $\begin{array}{c}
    \frac{1}{2} x^2
 \end{array}$
 &
 $\begin{array}{c}
    X^3_2 = x {\frac{\partial}{\partial{x}}}
            + u {\frac{\partial}{\partial{u}}}
            + 2\eta {\frac{\partial}{\partial\eta}},\\
    X^3_3 = e^{t}\left({\frac{\partial}{\partial{x}}}
            + {\frac{\partial}{\partial{u}}}
            - x {\frac{\partial}{\partial\eta}}\right),\\
    X^3_4 = e^{-t}\left({\frac{\partial}{\partial{x}}}
            - {\frac{\partial}{\partial{u}}}
            - x {\frac{\partial}{\partial\eta}}\right)
  \end{array}$
 &
 \begin{tabular*}{0.52\textwidth}{l|l}
    & \mmomconservtxt\\ \hline
    $X_1$ & \energconservtxt
    \\
    \hline
  $X^3_3$
  &
  $\begin{array}{c}
  D_x\big[
    e^t \big( u(x - \eta)(\eta + \frac{x^2}{2})
    - \frac{\eta}{2}(\eta + x^2)\big)
 \big]
 \\
 + D_t\big[
    e^t \big(x\eta - u(\eta + \frac{x^2}{2})\big)
 \big] = 0
  \end{array}$
  \\
  \hline
  $X^3_4$
  &
  $\begin{array}{c}
  D_x\big[
    e^{-t} \big( u(x + \eta)(\eta + \frac{x^2}{2}) + \frac{\eta}{2}(\eta + x^2)\big)
 \big]
 \\
 + D_t\big[
    e^{-t} \big(x\eta + u(\eta + \frac{x^2}{2})\big)
 \big] = 0
  \end{array}$
  \end{tabular*}
 \\
 \hline
 4 &
 $\begin{array}{c}
    -\frac{1}{2}x^2
 \end{array}$
 &
 $\begin{array}{c}
    X^4_2 = x {\frac{\partial}{\partial{x}}}
            + u {\frac{\partial}{\partial{u}}}
            + 2\eta {\frac{\partial}{\partial\eta}},\\
    X^4_3 =
                \cos{t}{\frac{\partial}{\partial{x}}}
                - \sin{t}{\frac{\partial}{\partial{u}}}
                + x\cos{t}{\frac{\partial}{\partial\eta}},\\
    X^4_4 = \sin{t}{\frac{\partial}{\partial{x}}}
                + \cos{t}{\frac{\partial}{\partial{u}}}
                + x \sin{t}{\frac{\partial}{\partial\eta}}
  \end{array}$
 &
 \begin{tabular*}{0.52\textwidth}{l|l}
    & \mmomconservtxt\\ \hline
    $X_1$ & \energconservtxt
    \\
    \hline
  $X^4_3$
  &
  $\begin{array}{c}
  D_t\big[
    u (\eta - \frac{x^2}{2}) \sin{t}
    - \eta x \cos{t}
 \big]
 \\
 + D_x\big[
    \left(
        u^2 (\eta - \frac{x^2}{2}) + \frac{\eta}{2} (\eta - x^2)
    \right) \sin{t}
    \\
    - xu(\eta - \frac{x^2}{2}) \cos{t}
 \big] = 0
  \end{array}$
  \\
  \hline
  $X^4_4$
  &
  $\begin{array}{c}
  D_t\big[
    -u (\eta - \frac{x^2}{2}) \cos{t}
    - \eta x \sin{t}
 \big]
 \\
 + D_x\big[
    -\left(
        u^2 (\eta - \frac{x^2}{2}) + \frac{\eta}{2} (\eta - x^2)
    \right) \cos{t}
    \\
    - xu(\eta - \frac{x^2}{2}) \sin{t}
 \big] = 0
  \end{array}$
  \end{tabular*}
 \\
 \hline
5 &
  $\begin{array}{c}
    x^{c},
    \\
    c\neq 0,1,2
 \end{array}$
 &
 $\begin{array}{r}
    X^5_2 = (2 - c) t {\frac{\partial}{\partial{t}}}
            + 2 x {\frac{\partial}{\partial{x}}}
            + c u {\frac{\partial}{\partial{u}}}
            + 2 c \eta {\frac{\partial}{\partial\eta}}
    \end{array}$
 &
    \multirow{2}{*}{
        Mass, momentum and energy ($X_1$) conservation
    }
 \\
\cline{1-3}
6 &
  $\begin{array}{c}
    e^x
  \end{array}$
  &
  $X^6_2 = t {\frac{\partial}{\partial{t}}}
            - 2{\frac{\partial}{\partial{x}}}
            - u {\frac{\partial}{\partial{u}}}
            - 2 \eta {\frac{\partial}{\partial\eta}}$
  &
  \\
\cline{1-3}
7 &
  $\begin{array}{c}
    c \ln |x|,
    \\ c \neq 0
  \end{array}$
  &
   $X^7_2 = t {\frac{\partial}{\partial{t}}}
            + x {\frac{\partial}{\partial{x}}}
            - c {\frac{\partial}{\partial\eta}}$
  &
  Mass, momentum and energy ($X_1$) conservation
  \\
\hline
\end{longtable}
\egroup

\begin{longtable}[c]{|l|c||c|c||c|c|}
\caption{Symmetries and conservation laws for finite-difference systems in Eulerian and Lagrangian coordinates}
\label{tab:HtableDelta}\\
\hline
\multirow{2}{*}{\#} & \multirow{2}{*}{Bottom $H$}
& \multicolumn{2}{c||}{Eulerian coordinates}
& \multicolumn{2}{c|}{Lagrangian coordinates}
\\ \cline{3-6}
                    &
                    & Symmetries
                    & Conservation laws
                    & Symmetries
                    & Conservation laws
\\ \hline\hline
\endfirsthead
\multicolumn{6}{c}%
{\emph{Table \thetable\ continued from previous page}} \\
\hline
\multirow{2}{*}{\#} & \multirow{2}{*}{Bottom $H$}
 & \multicolumn{2}{c||}{Eulerian coordinates}
  & \multicolumn{2}{c|}{Lagrangian coordinates}
  \\ \cline{3-6}
  &
  & Symmetries
  & Conservation laws
  & Symmetries
  & Conservation laws
\\ \hline\hline
\endhead
-- & $\begin{array}{c}\text{arbitrary}\\H(x)\end{array}$
& $X_1=\frac{\partial}{\partial t}$
& \multirow{3}{*}{
    \bgroup
    \def\arraystretch{1.25}
    $\begin{array}{l}
    \text{Scheme~(\ref{scm_arb}):}
    \\
    \empconservtxt
    \\
    \dtp{D}\left(
        u^2 (\eta + H) + \eta_+^2
    \right) +\\
     + \dhp{D}\Big[
        h\dtp{D}(\eta)u^2
        + \frac{1}{2}(\hat{\eta} + H) \times
        \\
        \times (u + \hat{u})(u^2 + \eta + \hat{\eta})
    \Big]=0
    \end{array}$
    \egroup
    }
    & $\frac{\partial}{\partial{t}}, \frac{\partial}{\partial{s}}$
    &
    \energyOrMomentumText
    \\
    \cline{1-3} \cline{5-6}
1
& $0$
& $X^1_2, X^1_3, \frac{\partial}{\partial{x}}$ & &
$\begin{array}{c}\frac{\partial}{\partial{x}}, t\frac{\partial}{\partial{x}},\\3t\frac{\partial}{\partial{t}}+2x\frac{\partial}{\partial{x}},\\3s\frac{\partial}{\partial{s}}+x\frac{\partial}{\partial{x}}\end{array}$
&
    \bgroup
    \def\arraystretch{1.25}
    \begin{tabular}[c]{@{}l@{}}
        Mass, center of mass, momentum\\
        and energy conservation ---\\
        see scheme~(\ref{schemeLagrMain})
    \end{tabular}
    \egroup
\\
\cline{1-3} \cline{5-6}
2 & $x$
& $X^2_2$
& &
\multicolumn{2}{c|}{Conservative scheme connected with scheme~(\ref{schemeLagrMain}) by (\ref{linearbtmTrLagr})}
\\ \hline
3 & $\frac{x^2}{2}$
& $X^3_2$
& \multirow{5}{*}{\noaddcls}
& $\begin{array}{c}X^3_2,\\e^t \frac{\partial}{\partial{x}},\\e^{-t}\frac{\partial}{\partial{x}}\end{array}$
& $\begin{array}{l}
    \text{Mass and energy conservation, and}
    \\
{\dtm{D}}(x\!\dtp{D}(e^t) - e^t {x}_t)
    -\dsm{D}\left(
        (\hat{x}_s \check{x}_s)^{-1} e^t
    \right)
    =0,
    \\
{\dtm{D}}(x\!\dtp{D}(e^{-t}) + e^{-t} {x}_t)
    -\dsm{D}\left(
        (\hat{x}_s \check{x}_s)^{-1} e^{-t}
    \right) = 0.
    \\
    \text{See scheme~(\ref{scheme_sq_lagr})}
   \end{array}$
\\
\cline{1-3} \cline{5-6}
4
&
$-\frac{x^2}{2}$
&
$X^4_2$
& &
$\begin{array}{c}
    X^4_2,\\
    \sin{t} \frac{\partial}{\partial{x}},\\
    \cos{t}\frac{\partial}{\partial{x}}
\end{array}$
&
$\begin{array}{l}
        \text{Mass and energy conservation, and}
        \\
        {\dtm{D}}(x\!\dtp{D}(\sin t)- {x}_t \sin t)
            -\dsm{D}\left(
                (\hat{x}_s \check{x}_s)^{-1}\sin t
            \right)
            = 0,
            \\
            \dtm{D}(x\!\dtp{D}(\cos t)- {x}_t \cos t)
            -\dsm{D}\left(
                (\hat{x}_s \check{x}_s)^{-1}\cos t
            \right) = 0.
            \\
            \text{See scheme~(\ref{scheme_sq_lagr2})}
  \end{array}$
\\
\cline{1-3} \cline{5-6}
5
& \begin{tabular}[c]{@{}l@{}}
    $\begin{array}{c}
        \\
        x^{c},\\
        c\neq 0,1,2
        \\
    \end{array}$
\end{tabular}
& $X^5_2$
&  \multirow{2}{*}{}
& $X^5_2$
&
    \bgroup
    \def\arraystretch{1.25}
    \begin{tabular}[c]{@{}l@{}}
        Mass and energy conservation.\\
        An additional conservation\\
        law may exist for $H=x^{-4/3}$\\
        only~(see Remark~\ref{rem:H-4/3})
    \end{tabular}
    \egroup
\\ \cline{1-3} \cline{5-6}
6
& $e^x$
& $X^6_2$
& &
\multicolumn{2}{c|}{\multirow{3}{*}{\energyOrMomentumText}}
\\ \cline{1-3}
7
& $\begin{array}{c}
     c\ln|x|,\\
     c\neq 0
   \end{array}$ & $X^7_2$
& & \multicolumn{2}{c|}{}
\\ \hline
\end{longtable}
\end{landscape}

\end{document}